\documentclass[12pt,letterpaper]{amsart}
\usepackage{amsmath, amsthm, latexsym, amsmath, amssymb, amscd, amsfonts, url, color, mathrsfs, mathtools, graphicx, xy, tikz-cd, float, fullpage
}

\newcommand{\op}{\operatorname}

\newcommand{\Z}{\mathbb{Z}}

\newcommand{\C}{\mathbb{C}}

\newcommand{\bbP}{\mathbb{P}}
\newcommand{\Th}{\operatorname{Th}}

\newcommand{\old}[1]{{}}

\newcommand{\bbO}{\mathcal{O}}
\newcommand{\degG}{\operatorname{deg}^{\text{G}}}

\newcommand{\unit}{\textbf{1}}
\newcommand{\Vect}{\op{Vect}}

\newcommand{\Tr}{\op{Tr}}
\newcommand{\Sp}{\text{Sp}}
\newcommand{\Top}{\text{Top}}
\newcommand{\res}{\op{res}}
\newcommand{\Sing}{\op{Sing}}
\newcommand{\Sym}{\op{Sym}}
\newcommand{\adj}[4]{#1\negmedspace: #2\rightleftarrows #3:\negmedspace #4}
\newcommand{\Sphere}{\mathbb{S}_G}
\newcommand{\D}{\mathbb{D}}
\newcommand{\Perf}{\text{Perf}}
\newcommand{\Jac}{\op{Jac}}
\newcommand{\Hom}{\op{Hom}}
\newcommand{\tr}{\op{tr}}

\newcommand{\CP}{\mathbb{C}\mathbb{P}}

\makeatletter
\newcommand{\tpitchfork}{%
  \vbox{
    \baselineskip\z@skip
    \lineskip-.52ex
    \lineskiplimit\maxdimen
    \m@th
    \ialign{##\crcr\hidewidth\smash{$-$}\hidewidth\crcr$\pitchfork$\crcr}
  }%
}
\makeatother

\DeclareMathOperator*{\colim}{\text{colim}}

\numberwithin{equation}{section}

\theoremstyle{remark}
\newtheorem{remark}[equation]{Remark}

\theoremstyle{plain}
\newtheorem{theorem}[equation]{Theorem}
\newtheorem{lemma}[equation]{Lemma}
\newtheorem{proposition}[equation]{Proposition}
\newtheorem{corollary}[equation]{Corollary}

\newtheorem{defn}[equation]{Definition}

\newtheorem{exa}[equation]{Example}

\definecolor{forestgreen}{rgb}{0.0, 0.27, 0.13}

\newcommand{\hidden}[1]{\footnote{Hidden:  #1}}
\renewcommand{\hidden}[1]{}

\begin{document}
\title{The equivariant degree and an enriched count of rational cubics}

\author[Bethea]{Candace Bethea}
\author[Wickelgren]{Kirsten Wickelgren}

\begin{abstract}
We define the equivariant degree and local degree of a proper $G$-equivariant map between smooth $G$-manifolds when $G$ is a compact Lie group 
    and prove a local to global result. We show the local degree can be used to compute the equivariant Euler characteristic of a smooth, compact $G$-manifold and the Euler number of a relatively oriented $G$-equivariant vector bundle when $G$ is finite. As an application, we give an equivariantly enriched count of rational plane cubics through a $G$-invariant set of 8 general points in $\CP^2$, valued in the representation ring and Burnside ring of a finite group. When $\Z/2$ acts by pointwise complex conjugation this recovers a count of real rational cubics. 
\end{abstract}
\maketitle 

\begingroup
\setlength{\parskip}{3pt}
\renewcommand{\baselinestretch}{0.25}\normalsize
\tableofcontents
\renewcommand{\baselinestretch}{1.0}\normalsize
\endgroup

\section{Introduction}\label{section:Introduction}

The topological degree of a continuous map $f\colon X\to Y$ of connected, oriented manifolds of the same dimension $n$ is the integer $d$ such that the induced map on top homology, $f_*\colon H_n(X)\to H_n(Y)$, is given by $f_*[X]=d[Y]$ in $H_n(Y)$. The degree is denoted $\deg^{\text{top}}(f)$, and if $X$ and $Y$ are $n$-spheres then $f_*\colon \Z\to \Z$ is multiplication by $d$. 
The local degree of $f$ at a point $x$ in $X$ in the preimage of a regular value $y$ in $Y$, denoted $\deg_x^{\text{top}}f$, is defined by first choosing an open set $U\subset X$ containing $x$ and $V\subset Y$ containing $y$ such that $U\cap f^{-1}(y)=\{x\}$ and $f(U)\subset V$. The local degree, $\deg_x^{top}f$, is  the degree of the induced map on $n$-spheres obtained by collapsing $\partial U$ and $\partial V$: 
$$ S^n
\simeq U/(U\setminus\{x\})\to V/(V\setminus \{y\})
\simeq S^n.$$
One of the triumphs of classical degree theory is that the degree of $f\colon X\to Y$ can be computed locally as a sum of local degrees. The local to global principle for the topological degree of a smooth map of smooth compact manifolds is due to Brouwer, stating
\begin{equation}\label{eq:non-G-brouwer-local}
\deg^{\text{top}}f =\sum_{x\in f^{-1}(y)}\deg_x^{\text{top}}f,
\end{equation} and this is independent of choice of regular value $y\in Y$. 
Our work gives definitions of a global and local equivariant degree that satisfy an analogous local to global principle for proper maps which are suitably oriented relative to an equivariant cohomology theory defined by an equivariant ring spectrum.

\begin{theorem}\label{theorem:Glocalglobal}
Let $E$ be a genuine equivariant $G$-spectrum with ring structure up to homotopy, and let $f: X \to Y$ be a proper $G$-equivariant map of smooth $G$-manifolds of the same dimension. Let $i_y\colon y\hookrightarrow Y$ be a regular value of $f$. Suppose that $f$ and $f|_x$ for $x \in f^{-1}(y)$ are equipped with compatible $E$-orientations.  
    Then there is an equality in $E^0_{G_y}(y)$,
    \begin{equation}
    i^*_y\deg^{G}f=\sum_{\substack{G_y\cdot x,\\ x\in f^{-1}(y)}}\tr_{G_x}^{G_y}\deg_x^{G_x}f.\end{equation}
    \end{theorem}

All terminology and notation are defined in Section \ref{section:conventions}. 

The topological degree appears in enumerative geometry in that it can be used to compute 
Euler numbers and Euler characteristics. A natural question is how to study enumerative results 
 under the presence of a group action, 
which one might pursue by replacing non-equivariant topological degrees with equivariant degrees. Quadratically enriched degree theory using $\mathbb{A}^1$-homotopy has been systematically developed by Morel, Kass-Wickelgren, and Kass-Levine-Wickelgren-Solomon in \cite{Morel-06, KW-local-degree, KLSW-rational-curves}, and has applications to enriched enumerative geometry more broadly. See for example \cite{27lines, Hoyois-lefschetz, SW-4-lines, larson-vogt-bitangents, quadratic-tropical, quadratic-twisted-cubics}. 
We develop equivariant degree theory for use in equivariantly enriched enumerative geometry and give an application to counting rational cubics in $\CP^2$. 

In stable equivariant homotopy, the Burnside ring $A(G)$ of a finite group $G$ replaces the integers in that it distinguishes homotopy classes of endomorphisms of representation spheres. $A(G)$ is defined to be the group completion of the monoid of isomorphism classes of finite $G$-sets with addition rule given by disjoint union, which also has a ring structure given by cartesian product. The Burnside ring is a free $\Z$-module with generators $[G/H]$ where $H$ runs over conjugacy classes of subgroups of $G$. When $G$ is a finite group, there is an existing definition of the equivariant degree of a self map of representation spheres $S^V\to S^V$, which defines an isomorphism $$\degG\colon \left[S^V, S^V\right]^G \xrightarrow{\sim} A(G).$$ This is analogous to $\op{deg}\colon \left[S^n,S^n\right]\xrightarrow{\sim} \Z$ being an isomorphism non-equivariantly, see \cite{s70}. 
The equivariant degree presented in Section \ref{section:degree} generalizes this by defining a degree for any proper $G$-equivariant map $f\colon X\to Y$ of smooth $G$-manifolds of the same dimension which is oriented relative to a genuine $G$-spectrum $E$. 
By computing $\deg^G f$ locally using the local equivariant degree constructed in Section \ref{section:localdegree}, we compute the degree as a sum of degrees of representation spheres, analogous to the non-equivariant local Brouwer degree formula. 

A brief overview of the construction of our degree and local degree runs as follows. Given $f\colon X\to Y$ a proper $G$-map of smooth $G$-manifolds of the same dimension and a genuine $G$-spectrum $E$ with ring structure up to homotopy, $f$ is \textit{relatively $E$-oriented} if there is an isomorphism \[\rho_X\colon \Th(-L_f)\wedge E\stackrel{\sim}{\to} X_+\wedge E,\] where $L_f\in K_0^G(X)$ is the tangent complex of $f$. We construct a proper pushforward in $E$-cohomology, $$f_*\colon E^{L_f}(X)\to E^0(Y).$$ 
When $f$ is relatively $E$-oriented, this induces an oriented pushforward 
\[f^{or}_*\colon E^0(X)\to E^0(Y).\]
We define the \textit{equivariant degree} of $f$ to be $\degG (f):= f^{or}_*(1)$ in $E^0(Y)$. 
Given a regular value $y$ in $Y$ of $f$ and $x\in f^{-1}(y)$, we obtain a map $f_x:=f|_x\colon x \to y$ which we assume is oriented and compatible with the relative $E$-orientation of $f$, which is the case when the orientation is parameterized (see Example  \ref{example:paramaterized_orientation_point}), thus $$\rho_x\colon \Th(-L_{fx})\wedge E\stackrel{\sim}{\to}x_+\wedge E.$$ 
Noting that $f_x$ is only $G_x:=\{g\in G\colon gx=x\}$ equivariant, this gives rise to another pushforward $$f_{x*}\colon E_{G_x}^{L_{fx}}(x)\to E_{G_x}^0(y)$$ that induces an oriented pushforward at $x$, $f^{or}_{x*}\colon E^0_{G_x}(x)\to E^0_{G_x}(y)$. The \textit{local degree of $f$ at $x$} is $\deg^{G_x}_x(f):=f^{or}_{x*}(1)$ in $E^0_{G_x}(y)$. 

These degree constructions are made explicit in Sections \ref{section:degree} and \ref{section:localdegree}, but the abbreviated constructions above make evident two important points. A first observation is that $f_x$ is only $G_x$-equivariant, 
so we need a suitable mechanism for obtaining a $G_y$-equivariant homotopy class from the $G_x$-equivariant class $f_{x*}^{or}(1)\in E^0_{G_x}(y)$. We use equivariant transfers, in particular $\tr_{G_x}^{G_y}\deg_x^{G_x}f$ encodes the $G_y$ orbit of $x\in f^{-1}(y)$ in so far as $\deg^{G_x}_xf$ contributes  to the global degree $\deg^Gf$. 

A second observation from the construction of the degree and local degree is that $f^{or}_*(1)$ and $f^{or}_{x*}(1)$ depend on $\rho_X 1\in E^{L_f}(X)$ and $\rho_x 1 \in E^{L_{f_x}}_{G_x}(x)$ respectively, i.e.,  $\degG_xf$ depends on a choice of relative $E$-orientation of $f$ in general. When $E$ is a complex oriented theory, this orientation is canonical.

One benefit of an equivariant local to global principle for the degree is that it can be used to compute the equivariant Euler number $n_G(V,u)$ of a relatively $E$-oriented equivariant vector bundle $V$ on a smooth $G$-manifold $X$ whose rank is equal to the dimension of $X$ with respect to a section $u\colon X\to V$. 
This is proved in Theorem \ref{thm:degree_computes_euler} in Section \ref{section:Euler}, restated below. 

\begin{theorem}\label{thm:poincare_hopf} 
   Let $X$ be a smooth, compact $G$-manifold of dimension $n$ and let $V$ be a rank $n$ $G$-vector bundle on $X$, relatively oriented with respect to a genuine equivariant ring spectrum $E$ for a finite group $G$. Let $H$ be a subgroup of $G$ and let $X\stackrel{u}{\to}V$ be an $H$-equivariant section with simple zeros. Let \[G/H\stackrel{u}{\longrightarrow}U\] be the corresponding map, with $U$ a finite dimensional space of sections with only simple zeros. There is an equality in $E^0(G/H)$
    \[
    n_H(V,u) = u^*\deg^G(f) = \sum_{\substack{H\cdot x,\\ u(x)=0}}\tr_{G_x}^{H}\deg_x^{G_x}(f).
    \]
    \end{theorem} 
The notation $f$ in the statement of Theorem \ref{thm:poincare_hopf} is the $E$-oriented map constructed from the universal section $\sigma\colon U\times X\to \pi_X^*V$, $\sigma(s,x)=s(x)$ by the composition $$\{\sigma = 0\}\to U\times X\to U.$$ 
See Section \ref{section:Euler} for a full treatment of the construction.

Brazelton uses parameterized homotopy to define the  equivariant Euler number $n_G(V,u)$ for equivariant vector bundles which are relatively $E$-oriented and of rank $n$ over smooth, compact (real or complex) $G$-manifolds of dimension $n$. When $E$ is complex oriented, Brazelton gives a formula for the equivariant Euler number as a sum of local indices. Theorem \ref{thm:poincare_hopf} gives a formula for local indices in terms of local degrees even when $E$ is not complex oriented. Furthermore, we prove an equivariant Gauss-Bonnet theorem using an unparameterized Euler class defined in Definition \ref{defn:euler_class_in_E}, compatible with a parameterized Euler class. 

\begin{theorem}Let $G$ be a finite group and let $X$ be a smooth, compact $G$-manifold. Let $\pi\colon X\to *$ be the map to a point. 
Then 
\[
\chi^G(X)=\pi_* e(TX)
\]in $\op{End}_{\Sp_G}(\Sphere)\cong A(G)$. 
\end{theorem}

Local degrees are valued in $\pi_0^G E\simeq E^0(*)$, and different genuine $G$-spectra $E$ enjoy various computational benefits. For example, when $E$ is the sphere spectrum $\Sphere$ and $G$ is a finite group, local degrees are valued in $\pi_0^G\Sphere\cong A(G)$, the Burnside ring, and we can take $H$-fixed points to recover orbit information for any $H\leq G$. When $E$ is complex equivariant $K$-theory $KU_G$, local degrees are valued in $\pi_0^G KU_G\cong R(G)$, the representation ring,  and we can use existing characteristic classes in equivariant $K$-theory to do enumerative computations. 

In Section \ref{section:example}, we make use of the results in Section \ref{section:Euler} to give an equivariant count of rational cubics through a $G$-invariant set $S$ of $8$ general points in $\CP^2$ for any finite group $G$, valued in either the Burnside ring $A(G)$ or the representation ring $R(G)$. The \emph{equivariant mass} $m^{G_C}(C)$ of a rational cubic $C$ is its $G_C$-equivariant Euler characteristic $\chi^{G_C}(C)$ where $G_C\leq G$ is the stabilizer of $C$. We show: 

\begin{theorem}\label{exampletheorem}
 Let $G$ be a finite group acting on $\CP^2$, and let $S$ be a $G$-invariant set of $8$ points in general position in $\CP^2$. 
    There exists a formula $N_{3, \CP^2, S}^G$ in $A(G)$ and $R(G)$ that counts orbits of rational plane cubics through $S$ in $\CP^2$.  
    Furthermore, \begin{equation}N_{3,\CP^2,S}^G=\sum_{\substack{ G\cdot C,\\ C \text{ a rational cubic }\\ \text{passing through }S}} \tr_{G_{C}}^G m^{G_{C}}(C), 
    \end{equation}
    where $G_C$ is the stabilizer of the rational cubic $C$. 
    \end{theorem}

When $\Z/2$ acts by pointwise complex conjugation, Theorem \ref{exampletheorem} recovers the signed count of real rational cubics through a set of $8$ general points that consists of real points and complex conjugate pairs of points, see Corollary \ref{cor:real_count}. In Section \ref{section:Hom_bundle} we show that an equivariant count of rational cubics can be obtained as a sum of local degrees as an equivariant Euler number computation. Let $X$ be the $G$-invariant general pencil of cubics in which the rational cubics through $S$ interpolate, detailed in Section \ref{section:example}. We have natural projections $\pi_1\colon X\to \CP^1\cong\bbP V_1$ and $\pi_2\colon X\to \CP^2\cong \bbP V_2$ for $V_1$ and $V_2$ rank 2 and 3 $G$-representations respectively. 

\begin{theorem}
With the assumptions in Theorem \ref{exampletheorem}, there is an $R(G)$-valued count of rational cubics through $S$ given by 
\begin{align*}
\pi_Xe^{K_G}(W) &= \sum_{\substack{G\cdot x,\\ d\pi_1(x)=0}} \tr_{G_x}^{G}\deg_x^{G_x}(d\pi_1) \\ 
&= \wedge^3V_2\otimes(\Sym^3 V_2^* - \Sym^9 V_2^* + 2\Sym^{12}V_2^* - \Sym^{15}V_2^*) \\ 
&- V_2\otimes \wedge^3V_2\otimes (\Sym^2V_2^* - \Sym^8V_2^* +2\Sym^{11} V_2^* - \Sym^{14} V_2^*) \\ 
&+ \bbO_{\op{Spec}\C}\\
&+ \det V_2\otimes \wedge^3 V_2\otimes (\Sym^6 V_2^*-\Sym^{12}V_2^* + 3\Sym^{18}V_2^* -2\Sym^{21}V_2^*)
\end{align*}
where $W:=\op{Hom}(\pi_1^*T^*\bbP V_1, T^*X) \to X$ and $d\pi_1$ is the section of $W$ determined by $\pi_1$.  
\end{theorem}
The results in this paper are related to existing equivariant enumerative results. In \cite{betheapencil}, the first named author of this work gives a count of orbits of nodal conics in a $G$-invariant general pencil of conics in $\C\bbP^2$ for finite groups not isomorphic to $\Z/2\times \Z/2$ or $D_8$ with linear action on $\C\bbP^2$. In \cite{brazeuler}, Brazelton gives a count of orbits of the 27 lines on a smooth cubic surface which is symmetric under the action of a finite group. In \cite{BB-bitangents}, Bethea-Brazelton give formulas for orbits of bitangents to any symmetric, non-hyperelliptic smooth plane quartic. Each result gives equivariant enumerative counts valued in the Burnside ring, though none use equivariant degree theory. 

\begin{remark}
The study of equivariant degree theory in general precedes this work. Given a representation of a finite group $G$ on $\mathbb{R}^n$ and $f\colon (\mathbb{R}^n,0)\to (\mathbb{R}^n,0)$ a finite, $G$-equivariant map, James Damon defines an equivariant degree as $G$-representation and relates it to the degree formula of Eisenbud and H. Levine \cite{damon-G-degree}. G\c{e}ba-Krawcewicz-Wu define the equivariant degree of a continuous equivariant map of the form $f\colon V\times\mathbb{R}^n\to V$ for a finite dimensional orthogonal representation $V$ of a finite group $G$, valued in $A(G)$ \cite{GKW-degree}.  Ize-Massab\'{o}and-Vignoli  define the $G$-degree of an equivariant map $\mathbb{R}^m\to \mathbb{R}^n$, valued in $A(G)$ \cite{IMV-degree}. Our work differs in terms of maps we allow, groups we allow, and rings the degree can take values in. We consider proper equivariant maps $f\colon X\to Y$ between smooth $G$-manifolds of the same dimension (real or complex) for all compact Lie groups $G$, and our degree is a pushforward in the cohomology of any genuine ring $G$-spectrum $E$. As special cases, we can choose $E$ to be the sphere spectrum (resp. complex equivariant $K$-theory) to obtain a degree valued in $A(G)$ (resp. $R(G)$).  
\end{remark}

\begin{remark}
    Six functor formalisms for cohomology theories defined by genuine $G$-spectra with ring structure up to homotopy are a rich and active area of study, extending beyond our use of a pushforward to define an equivariant degree. Brazelton  defines a pushforward in parameterized equivariant homotopy theory and uses it to define an equivariant Euler number of significant interest in our work \cite{brazeuler}. Cnossen has developed twisted ambidexterity in the parameterized $\infty$-category of genuine $G$-spectra, generalizing the Wirthm\"{u}ller isomorphsims in equivariant stable homotopy theory that we use in Sections \ref{section:degree} and \ref{section:localdegree} \cite{Cnossen-ambidexterity}. While we take a different perspective for our specific purposes of defining an equivariant degree, we mention these perspectives given the strength of their results. 
\end{remark}

\noindent{\bf Acknowledgements.} Candace Bethea was supported by the National Science Foundation award DMS-2402099 for the duration of this work. Kirsten Wickelgren was supported by NSF DMS-2103838 and NSF DMS-2405191. We thank Thomas Brazelton and Bastiaan Cnossen for helpful comments on this work. 

\section{Background}\label{section:conventions}

Throughout Sections \ref{section:degree} and \ref{section:localdegree}, we assume $G$ is a compact Lie group, and in Sections \ref{section:Euler},  \ref{section:example}, and \ref{section:Hom_bundle} we assume $G$ is a finite group. Given based $G$-spaces $X$ and $Y$, a map $f\colon X\to Y$ is $G$-equivariant if $g\cdot f(x)=f(g\cdot x)$ for all $g\in G$ and $x\in X$. At times we call $f$ a $G$-map for brevity. Let $\mathcal{U}_G$ be a complete $G$-universe in which any finite dimensional $G$-representation embeds. 
For any finite dimensional representation $V$ of $G$,  $S^V$ will denote the one-point compactification of $V$ with $G$ acting trivially on the point at $\infty$, which we will take to be the basepoint of $S^V$. Much of our conventions and notation follow \cite[Section 2]{H-H-R}. 

 $\Top_G$ will denote the  category of compactly generated, weak Hausdorff left $G$-spaces with morphisms given by $G$-equivariant maps and weak equivalences given by maps which induce isomorphisms on $\pi_*^G$. The smash product makes $\Top_G$ into a closed symmetric monoidal category, with unit object given by the zero sphere $S^0$ with trivial $G$-action. $\Top_G$ has an associated $\infty$-category by localizing by weak equivalences. We won't  differentiate notation here, but $\mathcal{S}_G$ is often used in the literature to denote the localization. By Elmendorf's Theorem \cite{elm83} (see also \cite{stephan16, nilpotence-MNN}), a model for $\Top_G$ is $\text{Fun}(\mathcal{O}_G^{op}, \Top)$, where $\mathcal{O}_G$ is the $\infty$-category associated to the full subcategory of $\Top_G$ of objects of the form $G/H$ for closed subgroups $H$ of $G$. The pointed category $\Top_{G*}$ is defined similarly. 

We work stably, with  $(\Sp_G,\wedge, \Sphere)$ denoting the symmetric monoidal stable $\infty$-category of genuine $G$-spectra with unit the equivariant sphere spectrum $\Sphere$. The unstable category $\Top_{G*}$ and the stable category $\Sp_G$ are related by the suspension functor 
\[
\Sigma^\infty \colon \Top_{G*}\to \Sp_G
\]
that takes a $G$-space $X$ to the spectrum whose value on an orthogonal $G$-representation $V$ is
\[
\Sigma^\infty X(V):=S^V\wedge X.
\]
There is a zero space functor $\Omega^\infty\colon \Sp_G\to \Top_{G*}$ given by $\Omega^{\infty}X:= X(0)$ and an adjunction
\[
\adj{\Sigma^\infty}{\Top_{G*}}{\Sp_G}{\Omega^\infty}.
\]
In general, given an inclusion of finite dimensional $G$-representations $V\subseteq W$ in the complete universe $\mathcal{U}_G$ we have the suspension and loop space functors defined by  $\Sigma^{W-V}X=S^{W-V}\wedge X$ and $\Omega^{W-V}(X)$ respectively. 
As a matter of convention, given a $G$-space $X$ we will typically exclude the $\Sigma^\infty$ notation and simply write $X_+$ to denote $\Sigma^\infty X$ in $\Sp_G$ when it is clear that we are working stably. 

We note that Sections \ref{section:example} and \ref{section:Hom_bundle} make use of the (complex) equivariant $K$-theory spectrum, $KU_G$. Recall that $\pi_0^G KU_G \cong R(G)$, the complex representation ring. We will never write $R_{\mathbb{C}}(G)$ to clarify we are working in the complex representation ring, rather simply $R(G)$ as this is the case throughout. 

For any subgroup $H$ of $G$, we have a restriction functor denoted by $\res^G_H\colon \Sp_G\to \Sp_H$.  We will sometimes use the notation $i^*_HX_+$ to denote the $H$-spectrum obtained from restricting the $G$-action on $X_+$ to $H$. The notation $i_H^*$ is typically used in the literature more generally to denote a pullback along $i\colon H\to G$ even when $i$ isn't the inclusion of a subgroup. Restriction admits a left adjoint given by $X_+\mapsto G_+\wedge_{H} X_+$ see \cite[2.2.3]{H-H-R}.

 Let $(\Vect_G(X), \oplus, \unit)$ denote the symmetric monoidal category of isomorphism classes of $G$-equivariant vector bundles on a $G$-space $X$. 
 Given a $G$-space $X$, $K_0^G(X)$ will denote the group completion of $\Vect_G(X)$. 
 %We may equivalently take $K_0^G(X)$ to be homotopy classes of chain complexes of $G$-vector bundles on $X$ under the additional relation that two quasi-isomorphism classes of chain complexes  are equivalent if and only if they're isomorphic after a direct sum of acyclic complexes. 
For every $V\in \Vect_G(X)$ we define the Thom spectrum of $V$, $\Th(V)$, to be the suspension spectrum of the Thom space $\mathbb{P}(V\oplus \unit)/\mathbb{P}V$, and this can be extended to associate a Thom spectrum to every element of $K_0^G(X)$. By abuse of notation, we won't write $\Sigma^{\infty}\Th(V)$ to indicate that we mean the Thom spectrum of $V$ as opposed to the underlying Thom space of $V$ as it is always clear from context.  See \cite{FLM-picard, units-of-ring-spectra,ando-blumberg-gepner-parameterized, Bachmann2017NormsIM} for a  treatment of invertible spectra $\op{Sph}_G$ and Thom spectra respectively, in the parameterized setting, which is more general than is needed in this work. 

An important point of clarification is that the unstable space underlying the Thom spectrum $\Th(V)$, $\mathbb{P}(V\oplus \unit)/\mathbb{P}V$, is not the same as the fiberwise Thom space $\Sigma_X^V$, which we will also denote $\Th_X(V)$ or $\Sigma^V_X S^0_X$, though it can be obtained from $\Sigma_X^V$ by identifying all points at infinity. We will also not notationally differentiate between the fiberwise Thom space and its suspension spectrum, as it will be clear throughout. 

 Given a genuine $G$-spectrum $E$ and $V\in K_0^G(X)$, we will denote $E$-cohomology groups by 
 \begin{equation} \label{eq:defn_of_cohomology}
 E^V(X) := [\Th(-V), E]^G. 
 \end{equation}
Note $E^0(X)=[X_+,E]^G$. 
Given a subgroup $H$ of $G$, we denote $[\Th(-V), E]^H$ by $E^V_H(X)$.

In the presence of a six functor formalism, \eqref{eq:defn_of_cohomology} is compatible with the definition using fiberwise Thom spectra and parameterized homotopy groups, 
 \[
 E^V_X = [\Th_X(-V),\pi^*E]^G,
 \]
 where $\pi\colon X\to *$ is the unique map to a point. We note that $\pi^*$ admits a left-adjoint $\pi_!$, constructed in \cite[Example 2.1.8]{may-sigurdsson-parameterized}. For $(Y,q,p)\in \text{Ret}_G^X$ a retractive space over $X$, meaning a $G$-space $Y$ and an equivariantly commutative diagram of the form 
\begin{center}
    \begin{tikzcd} 
    X \arrow["q"]{r}\arrow[bend right=20,swap, "id_X"]{rr}& Y\arrow["p"]{r} & X 
    \end{tikzcd}
\end{center}
The parameterized pushforward of spaces is $\pi_!(Y)=Y/q(X)$ \cite[Example 2.1.8]{may-sigurdsson-parameterized}. In particular, for any $G$-equivariant vector bundle $W\to X$ we have $\pi_!\Th_X(W) = \Th(W)$ as spaces. This extends to spectra by \cite[Lemma 4.3.1]{cary-parameterized}, \cite[Example 2.22]{brazeuler}. 
Observe 
 \begin{align} \label{parameterized_compatibility}
  E^V_X &= [\Th_X(-V),\pi^*E]^G \\
  &= [\pi_!\Th_X(-V),E]^G \notag\\
  &= [\Th(-V), E]^G, \notag
 \end{align}
 showing compatibility.

Given our motivation to generalize the degree of self maps of spheres, we finish the section with a brief overview of classical equivariant degree theory for representation spheres following \cite{s70}, which is well-known to homotopy theorists but may be new to enumerative geometers. Let $G$ be a finite group. From a $G$-equivariant map $f\colon V\to V$ of $G$-representations, we obtain a corresponding $G$-map of representation spheres, $S^V\to S^V$. For any $G$-space $X$ and $G$-representation $V$, there is a suspension $S^V\wedge X$. 
Equivariant homotopy classes of endomorphisms of basepoint-preserving $G$-equivariant maps of representation spheres are denoted $$[S^V, S^V]^G.$$ Given an inclusion 
$V\subseteq W$ of $G$-representations in $\mathcal{U}_G$, there is a map $$[S^V, S^V]^G\to [S^{W}, S^{W}]^G$$
defined by sending the class of $f\colon S^V\to S^V$ to $S^W\simeq S^V \wedge S^{W-V} \stackrel{f\wedge 1}{\longrightarrow} S^V \wedge S^{W-V} \simeq S^W$ where $W-V$ is the orthogonal complement of $V$ in $W$.
We can take the colimit over the poset of finite dimensional sub-representations of $\mathcal{U}_G$ under inclusion, $$\colim_{V}[S^V, S^V]^G.$$ 
The colimit is taken over maps $[S^V, S^V]^G\to [S^W, S^W]^G$, $f\mapsto f\wedge \op{id}_{S^{W-V}}$, for $V\subseteq W$ and $W-V$ the orthogonal complement, and this colimit is obtained at a finite stage \cite{s70}. 
For $V$ large and $G$ finite there is an isomorphism $$\deg^G\colon [S^V, S^V]^G\xrightarrow{\sim} A(G)$$ for which we have the following commutative diagram for any subgroup $H$ of $G$
\begin{equation}\label{diagram:segal_degree}
\begin{tikzcd}
\left[S^V, S^V\right]^G \arrow[r, "\degG"] \arrow[d]
& A(G) \arrow[d, "\epsilon_H" ] \\
\left[(S^V)^H, (S^V)^H\right] \arrow[r,  "\deg^{\text{top}}"]
&  \Z
\end{tikzcd}
\end{equation}
where $\epsilon_H\colon A(G)\to \Z$ is defined as $[S]\mapsto \left|S^H\right|$ and $\deg^{\text{top}}$ is the topological degree described in the introduction. The fact that \eqref{diagram:segal_degree} commutes means that the equivariant homotopy class of a $G$-map $f\colon S^V\to S^V$ is determined by $\deg^{\text{top}}(f^H)$ as $H$ runs through subgroups of $G$. As mentioned in the introduction, the Burnside ring $A(G)$ is the Grothendieck group of isomorphism classes of finite $G$-sets with addition given by disjoint union and multiplication by Cartesian product, which has generators $[G/H]$ as a free $\Z$-module as $H$ runs over conjugacy classes of subgroups of $G$. 

When $G$ is a finite group and $V$ is a complex $G$-representation, the equivariant degree of a map  $f\colon S^V\to S^V$ which is transverse to 0 is represented by the class of the finite $G$-set $[f^{-1}(0)]$ in $A(G)$.  
Since $f^{-1}(0)$ is a finite $G$-set whenever $f$ is transverse to 0, $f^{-1}(0)$ can be written as a disjoint union of orbits. Thus  $[f^{-1}(0)]$ can be written as a sum of transitive $G$-sets $$[f^{-1}(0)]=\sum_{(H_i)}n_i[G/H_i]$$ in $A(G)$ where $(H_i)$ are the conjugacy classes of subgroups of $G$ and $n_i\in \Z_{\geq 0}$ is the number of orbits of $f^{-1}(0)$ stabilized by $H_i$. This represents $\deg^G f$ in $A(G)$. 

\section{The Equivariant Degree}\label{section:degree}

Let $G$ be a compact Lie group and $E$ a genuine $G$-spectrum with ring structure up to homotopy. Let $X$ and $Y$ be smooth manifolds of the same dimension with proper $G$-actions.

\begin{defn}
    For an equivariant map $f: X \to Y$ of smooth $G$-manifolds, define the {\em tangent complex} $L_f$ in $K^G_0(X)$ to be the class represented by the complex $L_f$ 
    \[
TX \stackrel{Tf}{\longrightarrow} f^{*}TY
    \] in homological degrees $0$ and $-1$.
\end{defn}

\begin{defn}\label{defn:orientedmap}
    A map $f: X \to Y$ of smooth $G$-manifolds is {\em $E$-oriented} by an equivalence $\Th({-L_f}) \wedge E \simeq X_+ \wedge E$ in $\Sp_G$. 
\end{defn}

Equivariant orientation theory is a rich and active area of study, with key foundational and recent developments in \cite{may-orientations, CMW01, CW92, bhattacharya-zou}. 
 Definition \ref{defn:orientedmap} orients a map relative to a ring spectrum rather than orienting a specific vector bundle. Definition \ref{defn:orientedmap} is a specific case of the more general notion of a relative $E$-orientation determined by a trivialization of an $E$-line bundle given in \cite{units-of-ring-spectra}. 

By a proper map, we mean a map of topological spaces such that the inverse image of any compact subset is compact. 

\begin{proposition}\label{prop:any_map_smooth_compact_G_mflds_is_smoothable_proper}
A proper equivariant map $f: X \to Y$ of smooth $G$-manifolds factors as $$X \stackrel{i}{\hookrightarrow} W \stackrel{\pi}{\longrightarrow} Y$$ where $i$ is a closed $G$-embedding and $\pi$ is a smooth, proper family of $G$-manifolds.
\end{proposition}

\begin{proof}
Let $f: X \to Y$ be a proper $G$-map of smooth $G$-manifolds. By \cite[Theorem 6.1 and Section 7.4]{Mostow-equiv_embed} and its extension to the smooth case given in \cite[Introduction theorem]{Popper-eqivariant_embeddings}, there exists a $G$-representation $V$ and an equivariant map $X \to V$ which is a smooth embedding.  Since $f$ is proper, this embedding determines an equivariant smooth embedding $i': X \to S^V$. Let $i: X \to Y \times S^V$ be $i = (i',f)$.  
\end{proof}

In the terminology of \cite[Definition 2.39]{brazeuler}, proper maps between smooth $G$-manifolds are smoothable, proper.

We include a short construction of the pushforward here without reference to parametrized homotopy theory and 6-functor formalisms. See \cite[Section 4]{ando-blumberg-gepner-parameterized} for a general treatment of the pushforward in the parameterized setting. 

\begin{proposition}\label{prop:defn_f_*_oriented}
A proper map $f: X \to Y$ of smooth $G$-manifolds induces a pushforward $f_*: E^{L_f +  \zeta}(X)  \to E^{\zeta}(Y) $ on cohomology for any virtual $G$-representation $\zeta$.
\end{proposition}

\begin{proof}
Let $f: X \to Y$ be a proper $G$-map of smooth $G$-manifolds and let $i: X \to Y \times S^V$ be the equivariant smooth embedding constructed in the proof of Proposition~\ref{prop:any_map_smooth_compact_G_mflds_is_smoothable_proper}. By \cite[VI Theorem 2.2]{Bredon-transformation_groups}, there is a tubular neighborhood $U$ of $X$ in $Y \times S^V$ and a genuine equivalence $\Th({N_X (Y \times S^V)})\simeq U/\partial U$ where $N_X (Y \times S^V)$ denotes the normal bundle and $\partial U$ denotes the boundary of $U$. The Thom collapse map thus defines a map
\[
Y_+ \to  \Sigma^{-V}\wedge \Th({N_X (Y \times S^V)})
\] 
in $\Sp_G$. 
In $K_0^G(X)$, we have an equality $ N_X (Y \times S^V) = f^*TY + V - TX $ producing an equivalence $\Th({N_X (Y \times S^V)})\simeq \Sigma^V\wedge \Th(-L_f)$. We thus obtain a map
\[
Y_+ \stackrel{\D f}{\longrightarrow} \Th({-L_f})
\] which defines the claim pushforward after applying $[-,E] \circ \Sigma^{-\zeta}$. 

We see this construction to be independent of the choice of embedding $X \to V$ as follows. Given two such embeddings, we obtain an embedding $X \to V \times V'$ isotropic to the embedding $X \to V \to V \times V'$ where $V \to V \times V'$ sends $v$ to $(v,0)$. The resulting Thom collapse map is the $V'$ suspension of the Thom collapse corresponding to $X \to V$. 
\end{proof}

Let $E$ be a genuine equivariant $G$-spectrum,  with a ring structure up to homotopy. Let $f: X \to Y$ be a proper map of smooth $G$-manifolds of the same dimension equipped with a $E$-orientation
\[ 
\Th({-L_f})\wedge E\stackrel{\rho_X}{\stackrel{\sim}{\to}}X_+\wedge E.
\] 
The $E$-orientation of $f$ defines an isomorphism $E^{L_f}(X) \cong E^0(X)$, and we obtain an oriented pushforward $f_*^{or}: E^0(X) \to E^0(Y)$ from the pushforward $E^{L_f}(X)\to E^0(Y)$ constructed in Proposition \ref{prop:defn_f_*_oriented},
\begin{center}
    \begin{tikzcd} 
    E^{L_f}(X)\arrow["f_*"]{r}  & E^0(Y) \\
    E^0(X)\arrow[swap, "f^{or}_*"]{ur} \arrow [swap, "\rho_X"]{u}& 
    \end{tikzcd} 
\end{center}
Let $1$ in $E^{0}(X)$ denote the unit of the ring structure.

\begin{defn}\label{defn:global_degree}
    The {\em degree} of $f$ is defined to be $f_*^{or}( 1)$ in $E^0(Y)$, denoted $\deg^Gf$.\end{defn}
By definition, $f_*^{or}(1)=f_* (\rho_X1)$.

\section{The Local Equivariant Degree}\label{section:localdegree}

Before defining the local equivariant degree, we will establish a base change result for $E$-oriented maps and introduce equivariant transfers, which will be used in the proof of Theorem \ref{thm:G_local_global_with_proof}  (Theorem \ref{theorem:Glocalglobal} in the introduction).

\begin{lemma}\label{prop:pullback_cotangent_is_quasi_iso} Given a commutative diagram of smooth $G$-manifolds
\begin{equation}\label{diagram:general_base_change}
\begin{tikzcd}
W \arrow["\tilde{g}"]{r} \arrow["\tilde{f}"]{d}
& X \arrow["f"]{d} \\
Z \arrow["g"]{r} 
&  Y
\end{tikzcd}
\end{equation} 
there is a natural map $L_{\tilde{f}}\to \tilde{g}^*L_f$. If $f$ and $g$ are transverse and \ref{diagram:general_base_change} is a pullback diagram, then $L_{\tilde{f}}\stackrel{\sim}{\to} \tilde{g}^*L_f$ is a quasi-isomorphism. 
\end{lemma}
\begin{proof} Since  $\tilde{f}^*g^*TY=\tilde{g}^*f^*TY$, $(T\tilde{g}, \tilde{f}^*Tg)$ defines a map from $L_{\tilde{f}}\to \tilde{g}^*L_f$
\[
\begin{tikzcd}
L_{\tilde{f}}:= TW \arrow["T\tilde{f}"]{r} \arrow["T\tilde{g}"]{d}
& \tilde{f}^*TZ \arrow["\tilde{f}^*Tg"]{d} \\
\tilde{g}^*L_f:= \tilde{g}^*TX \arrow["\tilde{g}^*Tf"]{r} 
&  \tilde{g}^*f^*TY
\end{tikzcd}\]

If $f$ and $g$ are transverse and  \ref{diagram:general_base_change} is a pullback, then $W$ is a smooth manifold. The cone of $(T\tilde{g}, \tilde{f}^*Tg)$, $C_{(T\tilde{g}, \tilde{f}^*Tg)}$, is 
\begin{equation}\label{eq:cone_of_pullback}
0\to TW\to \tilde{g}^*TX\oplus \tilde{f}^*TZ\to \tilde{g}f^*TY\to 0
\end{equation}
in $K_0^G(W)$. Since $f$ and $g$ are transverse, 
\[ 
\tilde{g}^*TX\oplus \tilde{f}^*TZ\to \tilde{g}f^*TY\to 0
\]
is exact. The kernel of the first map is equivalent to $TW$, hence 
\[ 
0\to TW\to \tilde{g}^*TX\oplus \tilde{f}^*TZ
\]
is exact. Since \eqref{eq:cone_of_pullback} is exact, $C_{(T\tilde{g}, \tilde{f}^*Tg)}$ is null-homotopic as a complex in $\Perf(W, G)$, the category of complexes of $G$-equivariant vector bundles on $W$ that are quasi-isomorphic to bounded complexes of vector bundles on $W$. 
\end{proof}

\begin{proposition}\label{prop:base_change}
(Transverse Base Change) Let $f\colon X\to Y$ and $g\colon Z\to Y$ be proper, transverse $G$-maps of smooth $G$-manifolds, and let 
\begin{equation}\label{diagram:base_change_push_pull}
\begin{tikzcd}
W \arrow["\tilde{g}"]{r} \arrow["\tilde{f}"]{d}
& X \arrow["f"]{d} \\
Z \arrow["g"]{r} 
&  Y
\end{tikzcd}
\end{equation}
be the corresponding pullback diagram. Let $E$ be a genuine $G$-spectrum. Then 
\[
\tilde{f}_*\tilde{g}^* = g^*f_*.
\] are equal maps
\[
E^{L_f+\zeta}(X) \to E^{\zeta}(Z)
\] for any virtual $G$-representation $\zeta$.
 \end{proposition}

\begin{proof}
Let $V_X$ and $V_W$ be genuine $G$-representations such that $X\stackrel{i_X}{\hookrightarrow} V_X$ and $W\stackrel{i_W}{\hookrightarrow} V_W$ are smooth embeddings, and let $V:= V_X\times V_W$. $X$ embeds smoothly into $V$ by $(i_x,0)$, though we will write $X\stackrel{i_x}{\hookrightarrow}V$ for simplicity. There is a map $i_x\tilde{g}\colon W \to V$, which may not be an embedding in general. Since $f$ and $g$ are transverse maps and \eqref{diagram:base_change_push_pull}  is a pullback square, 
$X\stackrel{(f, i_X)}{\longrightarrow} Y\times V$ and $W\stackrel{(\tilde{f}, i_X\tilde{g})}{\longrightarrow} Z\times V$ are equivariant smooth embeddings. To see that $(\tilde{f}, i_X\tilde{g})$ is an embedding, note that it is an immersion by the exactness of \eqref{eq:cone_of_pullback}. $(\tilde{f}, i_X\tilde{g})$ is injective because the pullback $W$ of $f$ and $g$ injects into the product $Z \times X$ via $(\tilde{f},\tilde{g})$ and $i_X$ is an embedding. Moreover, $(\tilde{f}, i_X\tilde{g})$ is proper because $\tilde{f}$ and the diagonal of $Z \times V \to Z$ is proper. 
An injective, proper, immersion is an embedding. From these embeddings, we have 
$Y_+\stackrel{\D_f}{\longrightarrow} \Th({L_f})$ and $Z_+\stackrel{\D_{\tilde{f}}}{\longrightarrow} \Th({-L_{\tilde{f}}})$.

For any $V\in K_0^G(X)$, we have $\tilde{g}^*\in K_0^G(W)$ and a canonical injection $\tilde{g}^*V\to V$ given by projection, which induces a map on Thom spectra $p_V\colon \Th({\tilde{g}^*V})\to \Th(V)$. To show $\tilde{f}_*\tilde{g}^* = f^*g_*$, it suffices to show 
\begin{equation}\label{diagram:dual_and_pullback}
\begin{tikzcd}
\Th(-L_{\tilde{f}}) \arrow["p_{-L_f}"]{r} 
& \Th(-L_f)  \\
Z_+ \arrow["g"]{r} \arrow["\D_{\tilde{f}}"]{u}
&  Y_+\arrow["\D_{f}"]{u}
\end{tikzcd}
\end{equation}
commutes. 

Smashing \eqref{diagram:dual_and_pullback} with $S^V$ gives 
\begin{equation}
\begin{tikzcd}
\Th(N_W(Z\times V)) \arrow[r]
& \Th(N_X(Y\times V))  \\
S^V\wedge Z_+ \arrow[r] \arrow["S^V\wedge \D_{\tilde{f}}"]{u}
&  S^V\wedge Y_+\arrow["S^V\wedge \D_{f}"]{u}
\end{tikzcd}
\end{equation}
using the fact that $N_X (Y \times V) \simeq f^*TY + V - TX$ and $N_W (Z \times V) \simeq f^*TZ + V - TW $ in $K_0^G(X)$ and $K_0^G(W)$ respectively. From the embeddings $X\stackrel{(f, i_X)}{\longrightarrow} Y\times V$ and $W\stackrel{(\tilde{f}, i_X\tilde{g})}{\longrightarrow} Z\times V$, writing $(f, i_X)'$ for $(\tilde{f}, i_X\tilde{g})$ to simplify notation, we have the following commutative diagram
\begin{equation}
\begin{tikzcd}\label{diagram:embedding_diagram}
W \arrow[hookrightarrow, "(f\text{, }i_X)'"]{rr} \arrow["\tilde{g}"]{d}
& & Z\times V \arrow["(g\text{, }1)"]{d} \\
X \arrow[hookrightarrow, "(f\text{, }i_X)"]{rr} 
& &  Y\times V
\end{tikzcd}
\end{equation}
of spaces. Note $(f, i_X)$ and $(g,1)$ above are transverse. The maps $\tilde{g}$ and $(g,1)$ induce maps on tangent bundles, $T\tilde{g}\colon TW\to \tilde{g}^*TX$ and $T(g,1)\colon T(Z\times V) \to (g,1)^*T(Y\times V)$ over $W$ and $Z\times V$ respectively. Pulling $T(g,1)$ back along $(f, i_X)'$, we have over $W$ 
\[\begin{tikzcd}
TW \arrow[ "T(f\text{, }i_X)'"]{rr} \arrow["T\tilde{g}"]{dd}
& & (f,i_X)'^*T(Z\times V) \arrow["(f\text{, }i_X)'^*T(g\text{, }1)"]{dd}\arrow[r] &  \text{cofib }T(f,i_X)'
\arrow[dd]\\
\\
\tilde{g}^*TX \arrow["\tilde{g}^*T(f\text{, }i_X)"]{rr} 
& &  (f, i_X)'^*(g,1)^*T(Y\times V)\arrow[r] &  \text{cofib }\tilde{g}^*T(f,i_X)
\end{tikzcd}\]
Since $(f, i_X)$ and $(g,1)$ are transverse and \eqref{diagram:embedding_diagram} is a pullback, 
\[\text{cofib }T(f,i_X)'\simeq \text{cofib }\tilde{g}^*T(f,i_X). \]
Furthermore since $(f, i_X)'$ and $(f,i_X)$ are embeddings, $\text{cofib }T(f,i_X)' 
\simeq N_W(Z\times V)$
and 
$\text{cofib }\tilde{g}^*T(f,i_X) 
\simeq \tilde{g}^* N_X(Y\times V)$. 

We have thus produced an isomorphism 
\[N_W(Z\times V)\simeq \tilde{g}^*N_X(Y\times V)\]
and hence our desired map on Thom spectra
\[
\Th(N_W(Z\times V))\to \Th(N_X(Y\times V)),
\]
which necessarily implies that \eqref{diagram:dual_and_pullback} commutes. Desuspending by $S^V$ and taking $[-, E]$, we conclude that $\tilde{f}_*\tilde{g}^* = f^*g_*$ as the following commutes

\begin{equation}
\begin{tikzcd}\label{diag:final_base_change}
E^{L_{\tilde{f}}}(W) \arrow["\tilde{f}_*"]{d}
& E^{L_f}(X)\arrow["f_*"]{d} \arrow[swap, "\tilde{g}^*"]{l}\\
E^0(Z) 
&  E^0(Y) \arrow["g^*"]{l}
\end{tikzcd}
\end{equation}
\end{proof}

We now define the local degree. As always, let $E$ be a genuine equivariant $G$-spectrum with ring structure up to homotopy. Let $f: X \to Y$ be a relatively $E$-oriented map of smooth $G$-manifolds of the same dimension. Let $y\in Y$ be a regular value of $f$ and let $x$ be any point in $X$ in $f^{-1}(y)$.  The stabilizer subgroups of $x$ and $y$ will be denoted $G_x$ and $G_y$. 
Sard's Theorem guarantees we can  find a regular value $y\in Y$ of $f$. However, we can't guarantee there exists a regular value of $f$ that is fixed by the $G$-action on $Y$. For this reason, our equivariant local to global result in Theorem \ref{thm:G_local_global_with_proof} is an equality in $E^0_{G_y}(y)$, not necessarily $E^0_G(y)$. In the case where $f$ has a regular value that is a $G$-fixed point, equation \eqref{eq:local_to_global} is a $G$-equivariant equality.

In the proof of Proposition \ref{prop:defn_f_*_oriented} we constructed a map $\mathbb{D}f\colon Y_+\to \Th(-L_f)$ to define the pushforward $f_{*}\colon E^{L_f+\zeta}(X)\to E^{\zeta}(Y)$ after applying $[-,E] \circ \Sigma^{-\zeta}$.  The map $f\colon X\to Y$ is oriented by $\rho_X\colon \Th(-L_f)\wedge E \stackrel{\sim}{\to} X_+\wedge E$, and we defined the equivariant degree of $f$ to be $\degG f=f_{*}(\rho_X 1)=f^{or}_*(1)$ under the oriented pushforward $f^{or}_{*}\colon E^0(X)\to E^0(Y)$. We follow the same procedure to 
define the local degree as another oriented pushforward, defining a dual $\D f_x$ and pulling back the global orientation of $f$ to $f|_{x}$ for every $x\in f^{-1}(y)$. 

Let $\tilde{f} := f|_{f^{-1}(y)}$, and consider the $G_y$-equivariant diagram
\begin{equation}\label{diag:whole_preimage}
\begin{tikzcd}
f^{-1}(y) \arrow[hookrightarrow, ""]{r} \arrow[swap,"\tilde{f}"]{d}
& X \arrow["f"]{d} \\
y \arrow[hookrightarrow, "i_y"]{r} 
&  Y
\end{tikzcd}
\end{equation}
 Note this satisfies the hypotheses of Proposition \ref{prop:base_change} since $f^{-1}(y)$ contains only regular points. In particular, the following diagram commutes
\begin{equation}\label{diagram:D1}
\begin{tikzcd}
\Th(-L_{\tilde{f}}) \arrow[r] 
& \Th(-L_f)  \\
y_+ \arrow[hookrightarrow,"i_y"]{r} \arrow["\D\tilde{f}"]{u}
&  Y_+\arrow["\D{f}"]{u}
\end{tikzcd}
\end{equation}
Proposition \ref{prop:defn_f_*_oriented} shows that $\D f$ is independent of choice of embedding $X\hookrightarrow V$. Since $y$ is a regular value, $f^{-1}(y)$ is a finite collection of isolated points in $X$ and $\Th(-L_f)\simeq \bigvee_{x\in f^{-1}(y)} \Th(-L_{f,x})$. 
Diagram \eqref{diagram:D1} is thus equivalent to
\begin{equation}\label{diagram:D2}
\begin{tikzcd}
\bigvee_{x\in f^{-1}(y)} \Th(-L_{f,x}) \arrow[r] 
& \Th(-L_f)  \\
y_+ \arrow["i_y"]{r} \arrow["\D\tilde{f}"]{u}
&  Y_+\arrow["\D{f}"]{u}
\end{tikzcd}
\end{equation}

For any $x\in f^{-1}(y)$, there is \begin{equation}\label{eq:stable_dual_at_point}
\D f_x: y_+\to \Th(-L_{f,x})
\end{equation} in $\Sp_{G_x}$, obtained by applying $\res^{G_y}_{G_x}$ and projecting to the $\Th(-L_{f,x})$ summand in the wedge. 
In fact, $\D f_x$ is the dual of the map $f_x:= f|_x$, i.e., $f_x\colon x\to y$.

Applying $[-,E]$ to $\D f_x$, we obtain a pushforward in cohomology \begin{equation}\label{eq:f*nonorient}f_{x*}\colon E^{-L_{f_x}}_{G_x}(x)\to E^0_{G_x}(y).\end{equation} 

Suppose now that $f$ is relatively $E$ oriented by $\Th(-L_f)\wedge E\stackrel{\rho_f}{\simeq} X_+\wedge E$. Note $i_x^*L_f\stackrel{\sigma}{\simeq} L_{f_x}$. 
Suppose that $\rho_{f_{x}}$ is an $E$-orientation of $f_x$ \[
\Th(-L_{f_x})\wedge E\stackrel{\rho_{f_x}}{\simeq} x_+\wedge E
\] which is {\em compatible} with $\rho_f$ in the sense that the diagram 
\begin{equation}\label{cd:compatible_orientation_f_x}
\begin{tikzcd}
E^{L_f}(X)  \arrow[swap,"i_x^*"]{d}
& \arrow["\rho_f"]{l} E^0(X) \arrow["i_x^*"]{d} \\
E^{L_{f_x}}(x) 
&  \arrow["\rho_{f_x}"]{l} E^0(x)
\end{tikzcd}
\end{equation} commutes. 

\begin{exa}\label{example:paramaterized_orientation_point}
    An equivalence $\Th_X(-L_f) \wedge_X E \simeq X_+ \wedge_X E$ of parametrized spectra over $X$ induces a compatible orientation of $f_x$ by pullback.
\end{exa}

A compatible orientation $\rho_{f_x}$ induces an oriented pushforward 
\[f_{x*}^{or}\colon E^0_{G_x}(x)\to E^0_{G_x}(y).\]
Let $1\in E^0_{G_x}(x)$ denote the unit so that $\rho_{f_x} 1\in E^{L_{f_x}}_{G_x}(x)$. 

\begin{defn}\label{defn:local_degree} Given a regular value $y$ in $Y$ and $x\in f^{-1}(y)$ a simple zero of $f$, the {\em local  degree} of $f$ at $x$ is $\deg^G f_x = f_{x*}^{or}(1)$ in $E^0_{G_x}(y)$, denoted $\deg^{G_x}_xf$. \end{defn}

As was the case for Definition \ref{defn:global_degree}, $f_{x*}^{or}(1)=f_{x*}(\rho_x 1)$. 

 The  local to global principle in Theorem \ref{theorem:Glocalglobal}  relies on realizing $\D \tilde{f}$ in \eqref{diagram:D2} as a wedge of equivariant transfers, whence realizing the pullback of $\degG f$ along $i_y^*$ as a sum of transfers over the $G_y$ orbits of $f^{-1}(y)$.   Given any $\gamma\colon S^V\to E(V)$ in 
 $\pi_0^H E =  \colim_{V} [S^V, S^V\wedge E]^H$ we can assume $V = i^*_HW$ for some $G$-representation $W$, recalling from Section \ref{section:conventions} that $i_H^*W$ is the underlying $H$-representation of $W$, so that $\gamma$ is equivalent to $S^{i^*_HW}\to E({i^*_HW})$. 
 \begin{defn}\label{defn:transfer} \cite{Schwede}
The {\em transfer of $\gamma$ from $H$ to $G$, $\tr_H^G\gamma$,} is the composition
 \[
S^W\stackrel{t_H^G}{\longrightarrow} G_+\wedge_H S^{i^*_HW}\stackrel{G_+\wedge_H\gamma}{\longrightarrow}G_+\wedge_H E(i^*_HW)\stackrel{\epsilon}{\longrightarrow}E(W).
 \]
 in $\pi_0^GE$, where  $t_H^G$ is a Thom-Pontrjagin collapse map and $\epsilon$ is the counit of the adjunction between $G_+\wedge_H -$ and $\res^G_H$ given by the action map. 
 \end{defn}
 A detailed treatment of equivariant transfers can be found in Section 2 of \cite{Schwede}. See also \cite{Cnossen-ambidexterity} for a generalization of the Wirthm\"{u}ller isomorphisms. Given any subgroup $H$ of $G$ and $G$-space $X$,  $$G_+\wedge_H i^*_H X\cong (G/H)_+\wedge X.$$ 
 In particular, 
 \[
 G_+\wedge_H S^{i^*_HW}\cong (G/H)_+\wedge S^W\cong \bigvee_{G/H} S^{W}
 \]
 for any $G$-representation $W$. 
 
If  we have a $G$-manifold $X$ and  $G$-vector bundle $V\to X$, then for any $x$ in $X$,  
 \begin{equation} \label{eq:wedge_over_orbits}
 \Th(V|_x)\wedge_{G_x}G_+\simeq \bigvee_{x'\in G\cdot x}\Th(V|_{x'}).
 \end{equation} 
 This can be seen immediately from the equivariant bijection  of underlying spaces $$\varphi\colon \Th(V|_x)\times_{G_x}G\to \bigvee_{x'\in G\cdot x}\Th(V|_{x'})$$ given by $\varphi(v,g)=(gv, gx)$ and $\varphi(\infty, g)=(\infty, gx)$, which factors through the quotient $$\frac{\Th(V|_x)\times_{G_x} G}{\infty\times_{G_x} G}\simeq \Th(V|_x)\wedge_{G_x} G_+$$ and is basepoint preserving at $\infty$. 

 If $f\colon X\to Y$ is an equivariant map of smooth $G$-manifolds of the same dimension, $y\in Y$ is a regular value of $f$, and $x\in f^{-1}(y)$, then \begin{equation}\label{eq2}
 \bigvee_{G_y\cdot x}\Th(-L_{fx})\simeq \Th(-L_{fx})\wedge_{G_x}(G_y)_+.
 \end{equation}
 Indeed, by \cite[Theorem 6.1 and Section 7.4]{Mostow-equiv_embed} there exists a $G$-representation $V$ and an equivariant map $i\colon X \to V$ which is a smooth embedding.  We thus have $(i,f)\colon X\to V\times Y$, and writing $N_{(i,f)X}(Y\times V)$ as $N_X(Y\times V)$ for simplicity, \eqref{eq:wedge_over_orbits} applied to $N_X(Y\times V)$ becomes
\begin{equation}\label{eq3}
\bigvee_{x'\in G_y\cdot x} \Th(N_{X}(Y\times V)|_{x'})\simeq \Th(N_X(Y\times V)|_x)\wedge_{G_x}(G_y)_+.
\end{equation}
Since $N_X(Y\times V)\simeq -L_f+V$ in $K^G_0(X)$, smashing \eqref{eq3} with $S^{-V}$ gives \eqref{eq2}. 

 We now have the machinery necessary to prove introduction Theorem \ref{theorem:Glocalglobal},  stated below. 

\begin{theorem}\label{thm:G_local_global_with_proof}
    Let $E$ be a genuine equivariant $G$-spectrum with ring structure up to homotopy, and let $f: X \to Y$ be a proper $G$-equivariant map of smooth $G$-manifolds of the same dimension. Let $i_y\colon y\hookrightarrow Y$ be a regular value of $f$. Suppose that $f$ and $f_x$ for $x \in f^{-1}(y)$ are equipped with compatible $E$-orientations.  
    Then there is an equality in $E^0_{G_y}(y)$,
    \begin{equation}\label{eq:local_to_global}
    i^*_y\deg^{G}f=\sum_{\substack{G_y\cdot x,\\ x\in f^{-1}(y)}}\tr_{G_x}^{G_y}\deg_x^{G_x}f.\end{equation}
\end{theorem} 

\begin{proof} 
 
By writing $f^{-1}(y)$ as a wedge of points in each of it's $G_y$ orbits,  \eqref{diag:whole_preimage} becomes: 
\[
\begin{tikzcd}
\bigvee_{\substack{G_y\cdot x,\\ x\in f^{-1}(y)} }\left(\bigvee_{ x'\in G_y/G_x} x'_+\right) \arrow["\tilde{f}"]{dd}\arrow["i_{f^{-1}(y)}"]{r}
& X_+ \arrow["f"]{dd} \\
\\
y_+ \arrow["i_y"]{r}
&  Y_+
\end{tikzcd}
\] 
The compatible $E$-orientations of the $f_{x'}$'s determine a unique $E$-orientation of $\tilde{f}$ up to homotopy by \eqref{eq2}. This diagram fits into a larger commutative diagram: 
\[\begin{tikzcd}
x_+  \arrow["f_x"]{dd} \arrow[hookrightarrow,r] &  \substack{G_y\cdot x_+\text{ }\simeq \\ G_{y+}\wedge_{G_x}x_+} \arrow[hookrightarrow,r] \arrow["f_{G_yx}"]{dd}&\bigvee_{\substack{G_y\cdot x,\\ x\in f^{-1}(y)} }\left(\bigvee_{ x'\in G_y/G_x} x'_+\right) \arrow["\tilde{f}"]{dd}\arrow["i_{f^{-1}(y)}"]{r}
& X_+ \arrow["f"]{dd} \\
\\
y_+ \arrow[r]& y_+\arrow[r] & y_+ \arrow["i_y"]{r}
&  Y_+
\end{tikzcd}
\]
where all vertical arrows are compatibly relatively $E$-oriented.

By Proposition \ref{prop:base_change}, 
\[
i_y^*f_* = \tilde{f}_*  (i_{f^{-1}(y)})^*.
\]
We thus have
\begin{align*}
i_y^*\degG f =i_y^*f_*(\rho_X 1_X) &=\tilde{f}_*  (i_{f^{-1}(y)})^* (\rho_X 1_X)\\
&= \tilde{f}_*(\bigvee_{\substack{G_y\cdot x,\\ x\in f^{-1}(y)} } (\rho_{G_yx}1_{G_yx}))\\
&=\sum_{\substack{G_y\cdot x,\\ x\in f^{-1}(y)} } f_{G_yx*}(\rho_{G_yx}1_{G_yx}).
\end{align*}
We will show that \[
\tr_{G_x}^{G_y}f_{x*}(\rho_x 1_x) = 
f_{G_yx*}(\rho_{G_yx}1_{G_yx}), 
\] which will prove the theorem.
We will first show that the pushforward can be constructed unstably on Thom spaces, then explicitly describe both  $f_{x*}^{or}(1)$ and $f_{G_yx*}^{or}(1)$ as unstable maps. 

The orientation $\rho_x$ implies the existence of a Thom class $\mu_\rho\colon \Th(-L_{f_x})\to E$ in $E^{L_{fx}}(x)$ such that cupping with $\mu_\rho$ defines a Thom isomorphism $E^*(x)\to E^{*}(\Th(-L_{fx}))$ 
\cite[Theorem 6.1]{May-equiv_orientations}.  We can also choose a $G_y$ representation $W$ and smooth embedding $X\hookrightarrow W$ so that $\Th(N_x(Y\times W))\simeq \Sigma^{i^*_{G_x}W}\Th(-L_{f_x})$ as $G_x$-bundles. Using this identification, we obtain an unstable $G_x$-equivariant map of Thom spaces,
 \begin{equation}
 \begin{tikzcd}\label{eq:unstable_dual_to_point}
 \Th(N_x(Y\times W)) \simeq \Sigma^{i^*_{G_x}W}\Th(-L_{f_x})\arrow["\Sigma^{i^*_{G_x}W}\mu_{\rho}"]{rr} & &  \Sigma^{i^*_{G_x}W}E=E(i^*_{G_x}W).
 \end{tikzcd}
 \end{equation}
We will use the notation $D f_x$ to refer to the above map in contrast the stable dual $\D f_x$ of $f_x$ defined in equation \eqref{eq:stable_dual_at_point}. Note if $M$ is any other $G_x$-representation, we obtain an unstable map $\Sigma^{M+i^*_{G_x}W}\mu_{\rho}\colon S^M\wedge \Th(N_x(Y\times W))\to E(M+i^*_{G_x}W)$.

Let $V$ be any other $G_x$-representation. We can choose $M$ such that $V\hookrightarrow M+i^*_{G_x}W$. Any $\gamma\colon S^V\to E(V)$ in $E^0_{G_x}(x)$ is then equivalent to $S^{M+i^*_{G_x}W}\to E(M+i^*_{G_x}W)$, where we note that $S^V\wedge x_+\simeq S^V$, likewise for $S^{M+i^*_{G_x}W}$. We have defined $f^{or}_{x*}(\gamma)$ in $E^0_{G_x}(y)$ unstably as the following composition
\[\begin{tikzcd}
    S^{M+i^*_{G_x}W}\wedge y_+ \arrow[swap, "id_{S^M}\wedge \Sigma^{i^*_{G_x}W}D f_{x}"]{dd} \arrow["f^{or}_{x*}(\gamma)"]{rrrr}& & & & E(M+i^*_{G_x}W) \\
    \\
    S^M\wedge \Th(N_{x}(Y\times W))
    \arrow[swap, "\Sigma^{M+i^*_{G_x}W}\mu_{\rho}"]{rrrr}& & & &  E(M+i^*_{G_x}W) \arrow[swap, " mult\text{ }\circ \text{ }(\gamma\wedge id_E)(M+i^*_{G_x}W)"]{uu}
\end{tikzcd}
\]
where the right vertical map  
\[
E(V)\wedge x_+\stackrel{\gamma\wedge id_E}{\longrightarrow} E(V)\wedge E(V)\stackrel{mult}{\longrightarrow}E(V),
\]
 is obtained by $\gamma\wedge id_E$ composed with the ring multiplication of $E$ for any $G_x$-representation $V$, applied to $M+i^*_{G_x}W$. 

We describe \[f^{or}_{G_yx*}\colon E^0_{G_y}(G_{y+}\wedge_{G_x}x_+)\to E^0_{G_y}(y)\] unstably in an similar way. A general element of $E^0_{G_y}(G_{y+}\wedge_{G_x}x_+)$ is the class of a stable map $G_{y+}\wedge_{G_x}x_+\stackrel{\gamma}{\rightarrow} E$, whence as before an unstable map 
\[
\begin{tikzcd}G_{y+}\wedge_{G_x}x_+\wedge E(V)\arrow["\gamma\wedge id_E"]{r} & E(V)\wedge E(V) \arrow["mult"]{r}& E(V)
\end{tikzcd}
\]
obtained by smashing with $E$ and composing with the ring multiplication of $E$ for any $G_y$-representation $V$. Thus for $M+W$ in particular, $f^{or}_{G_yx*}(\gamma)$ is given by the composition: 
\begin{equation}
\begin{tikzcd}
    S^{M+W}\wedge y_+ \arrow[swap, "id_{S^M}\wedge \Sigma^WD f_{G_yx}"]{dd} \arrow[ "f^{or}_{G_yx*}(\gamma)"]{rrrr}& & & & E(M+W) \\
    \\
    \begin{tabular}{c}
    $S^M\wedge \Th(N_{G_yx}(Y\times W))\simeq$ \\ 
    $S^M\wedge (G_y/G_x)_+\wedge \Th(N_x(Y\times W))$ \end{tabular} 
    \arrow[swap, "\substack{id_{S^M} \wedge  id_{G_y/G_x}\wedge \mu_\rho \simeq \\ G_{y+}\wedge_{G_x}\mu_\rho}"]{rrrr}& & & & (G_y/G_x)_+\wedge E(M+W) \arrow[swap, "mult\text{ }\circ \text{ }(\gamma\wedge id_E)(M+W)"]{uu}
\end{tikzcd}
\end{equation}

When $\gamma$ is $1$, the unit of $E$, the composition $mult \circ (\gamma\wedge id_E)$ is equivalent to the action map $(G_y/G_x)_+\wedge E(M+W)\stackrel{act}{\longrightarrow} E(M+W)$. 
Thus $f^{or}_{G_yx*}(1)$ is the diagram above, where the right vertical map is simply
\[
(G_y/G_x)_+\wedge E(M+W)\stackrel{act}{\longrightarrow} E(M+W).
\]

By construction, $\Sigma^W D f_{G_yx} = t^{G_y}_{G_x}$ above, where $t^{G_y}_{G_x}$ is the notation for the Thom-Pontrjagin transfer defined in \cite{Schwede}. Thus $G_{y+}\wedge_{G_x}\mu_\rho = \Tr_{G_x}^{G_y}(\mu_\rho)$, where $\Tr^G_H$ is the external transfer defined in \cite{Schwede}. 
We conclude that \[f^{or}_{G_yx*}(1)=\tr_{G_x}^{G_y}(\mu_{\rho}).\]

To finish the proof, note that $f_x\colon x_+\to y_+$ is homotopic to the identity map $*_+\to *_+$ 
and $\rho_x 1$ is the orienting Thom class $\mu_\rho$. Therefore $\tr_{G_x}^{G_y}(\mu_{\rho}) = \tr_{G_x}^{G_y}(f_{x*}(\rho_x(1))$ and
\[f_{G_yx*}(\rho_{G_yx}1)=\tr_{G_x}^{G_y}(f_{x*}(\rho_x 1)),\]
as desired. 
\end{proof} 

 By choosing $E$ to be the sphere spectrum so that the transfer is that defined by the Burnside ring Mackey functor, we recover the global degree defined by Segal in \cite{s70} when $G$ is a finite group and $f\colon S^V\to S^V$ is a $G$-map of finite dimensional representation spheres. Specifically, 
 Definition \ref{defn:global_degree} agrees with the classical equivariant degree of \cite{s70}, and we can compute $\deg^G(f)$ by \begin{equation}
i_0^*\deg^{G}(f)=\sum_{\substack{G\cdot x,\\ x\in f^{-1}(0)}}\tr_{G_x}^{G}\deg_x^{G}(f)\end{equation} in $A(G)$, assuming without loss of generality that $0\stackrel{i_0}{\hookrightarrow}S^V$ is a regular value of $f$.  

\begin{corollary}\label{cor:fixed_pt} 
Let $f\colon X\to Y$ be as in Theorem \ref{thm:G_local_global_with_proof}. If $y\in Y$ is a regular value of $f$ which is also a $G$-fixed point, then equation \eqref{eq:local_to_global} is a $G$-equivariant equality. \end{corollary}
\begin{proof} This follows from the fact that $G_y=G$ 
when $y$ is a fixed point.
\end{proof}

\begin{proposition}
(Homotopy Invariance) Let $f,g\colon X\to Y$ be proper $G$-maps of smooth $G$-manifolds which are relatively $E$-oriented. 
Assume $f$ and $g$ are $G$-homotopic by a relatively $E$-oriented homotopy $H\colon X\times I\to Y$ such that the relative orientations of $f$ and $g$ are compatible with that of $H$. Then 
\[
\degG(f)=\degG(g)
\]
and in particular, for any $y\in Y$,
\[
i_y^*\degG(f)=i_y^*\degG(g)
\]
in $E^0_{G_y}(y)$.
\end{proposition}

\begin{proof} 
The homotopy $H$ induces an canonical isomorphism of vector bundles $f^*TY\cong g^*TY$, whence a path $L_f\to L_g$ in $K_0^G(X)$. There is a canonical equivalence 
\[
\Th(-L_f)\stackrel{H'}{\simeq} \Th(-L_g)
\]
in $\Sp_G$ after applying $\Th$. 
The orientations of $f$ and $g$, $\rho_f$ and $\rho_g$, which were assumed to be compatible with $H$, are compatible with $H'$. Thus the following commutes: 
\begin{center}
    \begin{tikzcd} 
    E^0(\Th(-L_f)) \arrow[phantom, swap, sloped, shift left,"\simeq"]{d} \arrow["\stackrel{\rho_f}{\simeq}"]{r} & E^0(X)\arrow[phantom, sloped, "="]{d}\\
    E^0(\Th(-L_g)) \arrow[ "\stackrel{\rho_g}{\simeq}"]{r} & E^0(X). 
    \end{tikzcd} 
\end{center}
Therefore the larger diagram commutes: 
\begin{center}
    \begin{tikzcd}
        E^0(\Th(-L_f))\arrow[bend left=20, "\simeq"]{rr}\arrow[swap, "f_*"]{dr}\arrow["\stackrel{\rho_f}{\simeq}"]{r} & E^0(X)\arrow[bend left=12, "g^{or}_*"]{d}\arrow[swap, bend right=12,"f^{or}_*"]{d} & E^0(\Th(-L_g))\arrow["g_*"]{dl}\arrow[swap,"\stackrel{\rho_g}{\simeq}"]{l}\\
         & E^0(Y) & \\
    \end{tikzcd}
\end{center}
immediately implying that $$\degG(f)=\degG(g)$$ in $E^0(Y)$. The statement then follows by pulling back along $i_y\colon y\hookrightarrow Y$ for any $y\in Y$. 
\end{proof}

\begin{exa}
    The map $f\colon \C\to \C$ defined by $z\mapsto z^3$ is invariant under $G=\Z/2$ acting by pointwise complex conjugation, and $\deg^G(f) = [G/G]+[G]$ in $A(G)$ and $ 1+\rho_G$ in $R(G)$ where $\rho_G$ is the regular representation. 
\end{exa}

\section{The Equivariant Euler Characteristic and Euler Number}\label{section:Euler}

In this section we relate the equivariant degree to two important invariants in equivariant homotopy theory, the equivariant Euler characteristic and Euler number. We asssume $G$ is a finite group throughout this section. 

We use \cite{May01} for our definition of the equivariant Euler characteristic, following the general definition of the categorical Euler characteristic of a dualizable object in a symmetric monoidal category. Let $X\in \Sp_G$ be dualizable (e.g., the suspension spectrum of a finite $G$-CW complex) with dual $\D X$ and let $f\colon X\to X$ be a $G$-equivariant endomorphism of $X$. The \textit{trace} of $f$, denoted $\tau^G(f)$, is the composition
\[
\Sphere \stackrel{\eta}{\longrightarrow}X\wedge \D X\stackrel{f\wedge 1}{\longrightarrow} X\wedge \D X\stackrel{\tau}{\longrightarrow} \D X\wedge X\stackrel{\epsilon}{\longrightarrow}\Sphere
\]
in $\pi_0^G\Sphere$, where $\eta$ and $\epsilon$ are coevaluation and evaluation maps. 
One can show as in the non-equivariant setting that $\D X\simeq \Th(-TX)$ in $\Sp_G$. 
A thorough treatment of equivariant duality for smooth, compact $G$-manifolds can be found in \cite[Theorem III.5.1]{LMS86}.

\begin{defn}\cite[V.1.1]{LMS86}\label{defn:euler_characteristic}
The {\em Euler characteristic} of $X$ is $\tau^G(\op{id}_X)$, denoted $\chi^G(X)$. 
\end{defn} 

 Importantly, the equivariant Euler characteristic of a smooth, compact $G$-manifold $X$ is related to non-equivariant Euler characteristic through the topological degree. For any subgroup $H\leq G$, there is a ring homomorphism $d_H\colon \pi_0^G\Sphere\to \Z$ defined by mapping a representative of a homotopy class $f\colon S^V\to S^V$ to $\deg^{\text{top}}f^H$, where $f^H\colon S^{V^H}\to S^{V^H}$ is the induced map on genuine $H$-fixed points of $f$. The equivariant and non-equivariant Euler characteristics are related by the by the fact that $d_H\chi^G(X)=\chi^{\text{top}}(X^H)$, and in particular 
$d_e\chi^G(X)=\chi^{\text{top}}(X)$. See \cite[Corollary V.1.8]{LMS86}. 

We can say more when we represent $\chi^G(X)\in \pi_0^G\Sphere$ using representations in which $X$ embeds smoothly when $X$ is a smooth, compact $G$-manifold, following \cite[III Section 5]{LMS86}. 
Equivariant Whitney embedding guarantees there exists a $G$-representation $V$ such that $X$ embeds smoothly into $V$ \cite{Mostow-equiv_embed}. The coevaluation is then the composition
\[
\eta: S^V \to \Th(N_X V) \to X_+ \wedge \Th(N_X V)
\] of the Pontryagin-Thom collapse map and the Thom diagonal, where $\Th(N_X V)$ here is the unstable Thom space rather than the Thom spectrum. Let $z: X \to N_X V$ denote the zero section. The normal bundle to the composition
\[
X \stackrel{\Delta}{\to} X \times X \stackrel{z \times 1}{\to} N_X V \times X
\] is computed to be 
\[
N((z \times 1)\circ \Delta)) \cong \Delta^* N(z \times 1) \oplus N(\Delta) \cong N_X V \oplus TX \cong V \times X.
\] Thus there is an induced Thom collapse map
\[
t: \Th(N_X V) \wedge X_+ \to X_+ \wedge S^V
\]
 The evaluation
\[
\epsilon: \Th(N_X V) \wedge X_+ \stackrel{t}{\to} X_+ \wedge S^V \to S^V
\] is the composition of the Pontryagin-Thom collapse $t$ and the map induced by $X_+ \to S^0$.

The composition \begin{equation}\label{eq:thom_def_of_euler}
S^V \stackrel{\eta}{\to} X_+ \wedge \Th(N_XV) \to \Th(N_XV) \wedge X_+  \stackrel{\epsilon}{\to} S^V
\end{equation} is a morphism in $\Top_G$ which represents the value of $\chi^G(X)\in \pi_0^G\Sphere$ on $V$. As stated previously, the non-equivariant degree of \eqref{eq:thom_def_of_euler} computes the non-equivariant topological Euler characteristic  $\chi^{\text{top}}(X)$, and moreover the equivariant degree (Definition \ref{defn:global_degree}) of \eqref{eq:thom_def_of_euler} computes the equivariant Euler characteristic of $X$. This is independent of $G$-representation $V$ for which $X$ embeds smoothly. 

We now turn to the equivariant Euler number, showing  local degrees can be used to compute the equivariant Euler number defined by Brazelton in \cite{brazeuler}. Let $X$ be a smooth $G$-manifold of dimension $n$ and let $V\to X$ be a rank $n$ $G$-vector bundle. Let $E$ be a genuine $G$-spectrum with a ring structure up to homotopy.

Let $\Gamma(V)$  denote the space of differentiable sections of $V$, and choose $W$ to be any finite dimensional vector subspace of sections which is mapped to itself under the action of $G$. Let $d$ be the dimension of $W$. Within $W$, let $U$ denote the space of sections of $W$ which only have simple zeroes, i.e., $X\stackrel{u}{\to}V$ is in $U$ if $\det \Jac u(p)\neq 0$ whenever $u(p)=0$. Here $\det \Jac u(p)$ denotes 
the section of the line bundle $\Hom(\det TX, \det V) \to X$ induced by the Jacobian determinant of the map $TX \to V$ induced by $Tu: TX \to TV$ at a point $p$ satisfying $u(p) = 0$. (At such points, the tangent space $TV$ to $V$ splits canonically as $V \oplus TX$).  

By projecting from $U\times X\stackrel{\pi_X}{\to}X$ and pulling back $V$, we have 
\[ \begin{tikzcd}
     \pi_X^*V\arrow[r]\arrow[d] & V \arrow[d] \\
     U\times X \arrow["\pi_X"]{r} & X
 \end{tikzcd}\]
Observe first that the bundle $\pi_X^*V\to U\times X$ has a tautological section $\sigma\colon U\times X \to \pi_X^*V$ given by $\sigma(u,p)=u(p)$.  We can also encode the vanishing behavior of sections in $U$ by defining \[Z:=\{ (u,p)\colon \sigma(u,p)=u(p)=0\}\subseteq U\times X.\] Combining this, we have the following commutative diagram: 
\[\begin{tikzcd}
    & \pi_X^*V\arrow[r] \arrow[d] & V\arrow[d] \\
    Z\arrow[hookrightarrow,"i_Z"]{r}\arrow[bend right, swap, "f:=\pi_U\circ i_Z"]{dr} & U\times X\arrow[bend right,swap, "\sigma"]{u}\arrow["\pi_U"]{d}\arrow["\pi_X"]{r} & X \\
    & U &
\end{tikzcd}\]
where we have defined $f$ to be the composition \[Z\stackrel{i_Z}{\longrightarrow}U\times X\stackrel{\pi_U}{\longrightarrow}U.\] We will use the degree of $f$ to compute the equivariant Euler number of $V$ with respect to a section $u$ in $U$. First we set up the necessary conditions for $\deg^Gf$ to be defined, in particular we show that $Z$ is a manifold and $f$ is relatively $E$-oriented.

\begin{lemma}
    $Z$ is a manifold of dimension $d$.
\end{lemma}

\begin{proof}
Let $(u,p)$ be a point of $Z$. $U$ is an open subset of $W$ because $\Jac$ is continuous, whence $U$ is a smooth manifold of dimension $d$. Thus $U \times X$ is a smooth manifold of dimension $d+n$. The Jacobian matrix of $\sigma$ evaluated at $(u,p)$ and restricted to the tangent space of $X$ has determinant $\Jac u(p)$. By construction of $U$, $\Jac u(p) \neq 0$. Thus $Z$ is a manifold of dimension $d+n-n$ by the inverse function theorem.
\end{proof}

By definition, $L_f$ is given by the complex \[TZ\stackrel{Tf}{\to}f^*TU \]
in degrees $0$ and $-1$ in $K_0^G(Z)$. We wish to exhibit an explicit equivalence $\Th(L_f)\wedge E\simeq Z\wedge E$, assuming that $V$ is relatively $E$-oriented. 
 
\begin{defn}\label{defn:thomas_rel_oriented_bdle}  \cite[Definition 4.15]{brazeuler} A rank n bundle $V\to X$ over a $G$-manifold of dimension $n$ is \em{relatively $E$-oriented} if there is an isomorphism 
\[
\Sigma^V_X\pi_X^*E \simeq \Sigma^{TX}_X\pi_X^*E.
\]
\end{defn}
If $V$ is relatively $E$-oriented, we show that $L_f$ is relatively $E$-oriented in the sense of Definition \ref{defn:orientedmap}, by exhibiting $L_f$ as a sum of $L_{i_Z}$ and $L_{\pi_U}$ in $K_0^G(Z)$: 

\begin{proposition}\label{prop:V_oriented_implies_cotangent_oriented}
    Let $X$ be a smooth, compact $G$-manifold of dimension $n$ and let $V\stackrel{\pi_V}{\longrightarrow} X$ be a rank $n$ $G$-vector bundle for a finite group $G$. A relative $E$-orientation of $V$ as in \cite[Definition 4.17]{brazeuler} induces a canonical $E$-orientation of $f$ as in Definition \ref{defn:orientedmap}. 
\end{proposition}
\begin{proof}
Because $Z$ is the zero locus of $\sigma$, we have the exact sequence
\[
0 \to TZ \to i_Z^* T(U \times X) \to i_Z^* \pi_X^*V \to 0
\] The tangent space of $U \times X$ is the direct sum $T(U \times X) \cong \pi_U^* TU \oplus \pi_X^* TX$, whence there is an exact sequence
\[
0 \to TZ \to f^* TU \oplus i_Z^* \pi_X^* TX \to i_Z^* \pi_X^* V  \to 0
\] This defines an exact triangle
\[
 L_f \to i_Z^* \pi_X^* TX \to i_Z^* \pi_X^* V \to 
\] on $Z$. This in turn defines an equivalence of $Z$-fiberwise Thom spaces
\[
\Sigma_Z^{L_f \oplus i_Z^* \pi_X^* V} S^0_Z \simeq \Sigma_Z^{i_Z^* \pi_X^* TX} S^0_Z
\]
Smashing with $E$ and composing with the pullback by $i_Z^*$ of the equivalence $\Sigma_{X}^{TX} E \simeq \Sigma_{X}^{V} E$ produces an equivalence
\[
\Sigma_Z^{i_Z^* \pi_X^* V \oplus L_f} E \simeq \Sigma_Z^{i_Z^* \pi_X^* V} E
\]
Applying $\Sigma_Z^{-i_Z^* \pi_X^* V}$ defines an equivalence $\Sigma_Z^{L_f} E \simeq \Sigma_Z^{0} E$ which in turn determines a canonical $E$-orientation $\Th(L_f) \wedge E \simeq E$ of $f$ by passing to quotients from the fiberwise Thom spaces over $Z$ to Thom spaces.
\end{proof}

The equivariant Euler number of \cite[Definition 5.3]{brazeuler} is defined for $G$-equivariant vector bundles of rank $n$ over smooth, compact (real or complex) $G$-manifolds of dimension $n$ as a pushforward in the cohomology of any genuine $G$-ring spectrum $E$ so long as $V$ is relatively $E$-oriented. When $E$ is complex oriented and $V\to X$ is a complex vector bundle on a complex manifold, Brazelton gives a formula for the equivariant Euler number as a sum 
\[
n_G(V,s) = \sum_{G\cdot x\subseteq Z(s)} \tr_{G_x}^G(1)
\]
in $E^0(*)$, which is independent of choice of $G$-equivariant section $s\colon X\to V$, see \cite[Lemma 5.21]{brazeuler}. While the Euler number is defined for any genuine $G$-spectrum $E$, this formula for local indices in terms of transfers only holds when $E$ is complex oriented. We show that there is a formula for the local indices appearing in the definition of the equivariant Euler number in terms of local degrees even when $E$ is not complex oriented. 

\begin{theorem}\label{thm:degree_computes_euler}
Let $X$ be a smooth, compact $G$-manifold of dimension $n$ and let $V$ be a rank $n$ $G$-vector bundle on $X$, relatively oriented with respect to a genuine equivariant ring spectra $E$ for a finite group $G$. Let $H$ be a subgroup of $G$ and let $X\stackrel{u}{\to}V$ be an $H$-equivariant section with simple zeros. Let \[G/H\stackrel{u}{\longrightarrow}U\] be the corresponding map, where $U$ is a finite dimensional space of sections with only simple zeros. Then we have an equality in $E^0(G/H)$
    \begin{equation}\label{eq:degree_euler}
    n_H(V,u) = u^*\deg^G(f) = \sum_{\substack{H\cdot x,\\ u(x)=0}}\tr_{G_x}^{H}\deg_x^{G_x}(f).
    \end{equation}
\end{theorem}

\begin{proof}
   Let $\rho$ denote the relative orientation of $V$ and $f_{\rho}$ the induced orientation of $f: Z \to U$ of Proposition~\ref{prop:V_oriented_implies_cotangent_oriented}. By Theorem~\ref{thm:G_local_global_with_proof}, 
   \[
   u^*\deg^G(f) = \sum_{\substack{H\cdot x,\\ u(x)=0}}\tr_{G_x}^{H}\deg_x^{G_x}(f).
   \] The Euler number $n_H(V,u)$ is likewise a sum over the zeros of $u$,
   \[
   n_H(V,u) = \sum_{\substack{H\cdot x,\\ u(x)=0}} \op{ind}_x u.
   \] By assumption, $x$ is a simple zero of $u$, whence it follows that $\op{ind}_x u$ is $\tr_{G_x}^{H}$ applied to the $G_x$-degree of the map induced by the derivative of $u$
   \begin{equation}\label{eq:ind_e}
   S^{T_xX} \wedge E \to S^{V_x} \wedge E \stackrel{\rho}{\simeq} S^{T_xX} \wedge E 
   \end{equation} and the relative orientation $\rho$ of $V$. See
   \cite[Proposition 5.6]{brazeuler}. 

   The derivative of $f$ at $x$ is likewise full rank, whence $\deg_x^{G_x}(f)$ is the $G_x$-degree of the map induced by the derivative of $f$
   \begin{equation}\label{eq:locdegffore}
   S^{T_{(u,x)}Z} \wedge E \to S^{T_u U} \wedge E \stackrel{\rho_f}{\simeq}  S^{T_{(u,x)}Z} \wedge E 
   \end{equation} As in the proof of Proposition~\ref{prop:V_oriented_implies_cotangent_oriented}, we have a canonical isomorphism 
   \[
   T_{(u,x)}Z \oplus V_x \cong T_u U \oplus T_x X 
   \] Smashing \eqref{eq:locdegffore} with $S^{V_x - T_u U}$, we have that the $G_x$-degree of \eqref{eq:locdegffore} equals that of
   \[
   S^{T_xX} \wedge E \to S^{V_x} \wedge E \stackrel{\rho_f}{\simeq} S^{T_xX} \wedge E .
   \] This degree equals the $G_x$-degree of \eqref{eq:ind_e} by construction of $\rho_f$. Thus $\op{ind}_x u = \tr_{G_x}^{H}\deg_x^{G_x}(f)$ showing the claimed equality.
\end{proof}

Theorem \ref{theorem:Glocalglobal} thus implies that the equivariant Euler number can be computed as a sum of local degrees along the vanishing locus of an equivariant section. Note that when $u$ is a $G$-equivariant section of $V$, $u^*\deg^G(f)\in E^0_G(u)$. When $E$ is a complex oriented cohomology theory, $n_G(V,u)$ is independent of section $u$ by \cite[Lemma 5.4]{brazeuler}, whence so is the pullback of $\deg^G(f)$ along any $G$-equivariant section of $V$. 

We finish by relating the equivariant Euler number of the tangent bundle of a smooth, compact $G$-manifold $X$ to $\chi^G(X)$ when $G$ is a finite group and $TX$ is relatively oriented. As $\chi^G(X)$ is not defined using parameterized spectra and the Euler number is, we first propose a definition of an unparameterized Euler class compatible with the parameterized definition by \eqref{parameterized_compatibility}.

The parameterized Euler class of a $G$-equivariant vector bundle $W\to X$ can be represented by
$$e_X(W): 
X_+ \simeq \Th_X(0) \to \Th_X(W) = \Sigma^W_XS^0_X,$$ the map on fiberwise Thom spectra canonically induced by the zero section $s\colon X\to W$.  Given another vector bundle $W'\to X$, we can shift this class by $W'$
\[
\Sigma^{W'}_Xe_W\colon \Sigma^{W'}_X S^0_X\to\Sigma^{W+W'}_X S^0_X.
\]

Taking $W' = -W$, this gives a map of parameterized spectra
\[
\Th_X(-W) \to \Th_X(0)\simeq S^0_X. 
\]

 Let $\pi\colon X\to *$ be the map to a point. 

\begin{defn}\label{defn:unparameterized_euler_class}
The \emph{Euler class} of $W\to X$, denoted $e(W)$, is the map 
\begin{center}
    \begin{tikzcd}
        \Th(-W)\arrow["\pi_! \Sigma_X^{-W}e_X(W)"]{rrr} &&& \Th(0)\simeq X_+
    \end{tikzcd}
\end{center}
in $\Sp_G$. 
\end{defn}

We can also define the Euler class with respect to a ring spectrum. 

\begin{defn}\label{defn:euler_class_in_E} 
For any genuine $G$ spectrum $E$, there is a projection 
\[
X_+ \wedge E \to S^0\wedge E \to E.
\]
We define the {\emph{Euler class of $W\to X$ with respect to $E$}}, denoted $e^E(W)$, to be the composition 
\begin{equation}\label{eq:shifted_euler_class}
\Th(-W)\wedge E \stackrel{e(W)\wedge E}{\longrightarrow} X_+\wedge E \to E
\end{equation} 
in $\Sp_G$. 
\end{defn}

Finally, recall that projection to a point $X\stackrel{\pi}{\to}*$ induces an unparameterized pushforward $E^V(X)\stackrel{\pi_*}{\to} E^0(*)$, whose construction is outlined in detail in Section \ref{section:degree}. When the target of a map is a point, the pushforward in $E$-cohomology can intuitively be thought of as equivariant fiber integration. Composing $\pi_*$ with the Euler class in Definition \ref{defn:euler_class_in_E} defines an unparameterized \emph{Euler number} in $E^0(*)$, $\pi_* e^E(W)$.

\begin{theorem}\label{thm:gauss_bonet} Let $G$ be a finite group and let $X$ be a smooth, compact $G$-manifold. Let $\pi\colon X\to *$ be the map to a point. 
Then 
\[
\chi^G(X)=\pi_* e(TX)
\]in $\op{End}_{\Sp_G}(\Sphere)\cong A(G)$. 
\end{theorem} 
\begin{proof} 
  $\chi^G(X)$ is defined by \eqref{eq:thom_def_of_euler}. This is represented by the composition in $\Top_G$
  \[
  S^V \stackrel{t}{\to} \Th(N_X V) \to X_+ \wedge \Th(N_X V) \to \Th(N_X V) \wedge X_+ \stackrel{t}{\to} X_+ \wedge S^V \to S^V
  \] for $V$ a $G$-representation such that $X\hookrightarrow V$ is a smooth embedding, where again this is a map of underlying Thom spaces rather than Thom spectra. The composition
  \begin{equation}\label{Euler_piece_epsilon_eta}
  \Th(N_X V)\to X_+ \wedge \Th(N_X V) \to \Th(N_X V) \wedge X_+ \stackrel{t}{\to} X_+ \wedge S^V \simeq \Th(N_X V \oplus TX)
  \end{equation} admits the following alternate description. Let $s:X \to W $ be the zero section of a $G$-vector bundle $W \to X$. Then 
  \[
  s \times 1: X \to W \times X
  \] is an embedding with normal bundle
  \[
  N(s \times 1) \cong \Delta^*(W \boxplus 0) + N(\Delta) \cong W \oplus TX,
  \] where $\Delta: X \to X \times X$ denotes the diagonal. The embedding $s \times 1$ determines a Pontryagin--Thom collapse map 
  \[
  t(s): \Th(W) \wedge X_+  \simeq \Th(W \boxplus 0) \to \Th(N(s \times 1)) \simeq \Th(W \oplus TX).
  \] The composition 
  \[
  \Th(W) \to \Th(W) \wedge X_+ \stackrel{t(s)}{\to} \Th(W \oplus TX) 
  \] of the Thom diagonal with $t(s)$ is homotopic to $\Th(W \to W \oplus TX)$ the map on Thom spaces induced by the canonical injection $W \to W \oplus TX$. Applying this when $s$ is the zero section $s\colon X \to N_X V$ computes the map in \eqref{Euler_piece_epsilon_eta} to be homotopic to $\Th(N_X V \to N_X V \oplus TX)$. 

  In total, we thus have that $\chi^G(X)$ is represented by
  \begin{equation}\label{eq:normal_bdle_euler_characteristic}
  S^V \stackrel{t}{\to} \Th(N_X V) \to \Th(N_X V + TX)\simeq X_+ \wedge S^V \to S^V
  \end{equation} 
  where $\Th(N_X V)  \to  X_+ \wedge S^V$ is the map $\Th(N_X V \to N_X V \oplus TX)$ of Thom spaces. 
  Thus 
  \begin{equation}\label{eq:shifted_euler_class_in_proof}
  \Th(N_X V) \to \Th(N_X V + TX)\simeq X_+ \wedge S^V \to S^V
  \end{equation}
  is homotopic to $\Sigma^V e(TX)$ after taking the suspension of \eqref{eq:shifted_euler_class_in_proof}. Therefore \eqref{eq:normal_bdle_euler_characteristic} is indeed $\pi_*\Sigma^V e(TX) \simeq \pi_* e(TX)$. 
\end{proof} 

Since the equivariant degree of \eqref{eq:normal_bdle_euler_characteristic} (or equivalently \eqref{eq:thom_def_of_euler}) computes $\chi^G(X)$, it is immediate that $\pi_* e(TX)$ can be computed as a sum of local equivariant degrees by Theorem \ref{thm:G_local_global_with_proof}. We note that Theorem \ref{thm:gauss_bonet} should also have a more general formulation in terms of parameterized spectra.

\section{An Equivariant Count of Rational Cubics Through 8 Points}\label{section:example}

A degree $d$ rational plane curve is a map $u\colon \CP^1\to \CP^2$ defined by homogenous degree $d$ polynomials $f_0, f_1, f_2$ so that $[s,t]\stackrel{u}{\mapsto} [f_0(s,t), f_1(s,t), f_2(s,t)]$. There are infinitely many rational plane curves of a given degree in general. However we can impose the condition that $\text{im}(u)$ passes through a given set of marked points in general position in $\CP^2$, and there are finitely many rational degree $d$ plane curves passing through $3d-1$ general points. A groundbreaking success of Gromov-Witten theory was the formulation of a recursive formula in $d$ for the count of rational curves through $3d-1$ general points for $d\geq 5$  \cite{KM-recursive}. 

The first nontrivial case is that for cubic curves, which states that there are 12 rational cubics passing through any 8  points in general position in $\CP^2$. This case can be seen by first showing that the rational cubics through any 8 points are in bijection with the nodal cubics in a general pencil of plane cubics. The number of nodal cubics in a general pencil is then given by the degree of the top Chern class of the bundle of principle parts on $\mathcal{O}_{\CP^2}(3)$, with \[
\deg c_2(\mathcal{P}^1(3))=3(2)^2=12.
\]
See \cite[Chapter 7]{EH16} for a detailed proof of the count of nodal curves in a general degree $d$ pencil of curves in $\CP^n$.  

We explore this question under the presence of a finite group action. Let $S$ be a $G$-invariant set of 8 points in general position in $\CP^2$. Instead of asking for the count of rational cubics through $S$, we will ask for an equivariant count of orbits of rational cubics through $S$. We prove that there is a formula in both $A(G)$ and $R(G)$ 
of the orbits of rational cubics through $S$.  
By taking the cardinality (resp. rank) of the count in $A(G)$ (resp. $R(G)$) we obtain the integer 12, recovering the classical count. By taking the cardinality of $H$-fixed points of $N_{3,\CP^2, S}^G$, considered in $A(G)$, for any subgroup $H$ of the whole group, we recover the number of rational cubics stabilized by $H$. 

When $G=\Z/2$ acts on $\CP^2$ by pointwise complex conjugation, we can recover the signed integer count of real rational curves through a $G$-invariant set of 8 points from $N_{3,\CP^2, S}^G$.  When $N_{3,\CP^2, S}^G$ is considered as an element of $A(G)$,  $-|(N_{3,\CP^2, S}^G)^G|$ recovers the count of real rational cubics through $S$. When $N_{3,\CP^2, S}^G$ is considered as an element of $R(G)$, minus the character of the non-identity element $g\in G$, $-\chi_{N_{3, \CP^2, \Sigma}^G}(g)$, equals the number of real rational cubics through $S$. In orther words, by specificying a specific group and an action we recover both a complex and real rational curve count from the same formula, $N_{3, \CP^2, S}^G$. This is Corollary \ref{cor:real_count}. 

The mass of a nodal curve will encode the equivariant topology of the branches of a rational cubic, including the topology of the branches of a real rational cubic in the case when $\Z/2$ acts by complex conjugation. 
This generalizes the $A(G)$-valued weight of a node from  \cite{betheapencil}.

\begin{defn}\label{defn:mass}
The equivariant {\em mass} of a nodal curve $C$ is its equivariant Euler characteristic, \[m^{G_C}(C)= \chi^{G_C}(C),\] where $G_C$ is the stabilizer of $C$.
\end{defn}

When $E$ is the sphere spectrum and $M$ is any finite $G$-CW complex, the definition of $\chi^G(M)$ we use is that of Definition \ref{defn:euler_characteristic} in terms of the trace of $\text{id}_M$. When $E$ is $KU_G$ we denote by $\chi^{K_G}(M)$ the equivariant Euler characteristic in $K$-theory, which is the composition of $\chi^G(M) \in\text{End}_{\Sp_G}(\Sphere)$ and $\pi_0^G\Sphere\to \pi_0^GKU_G$. Since $\pi_0^G KU_G\cong R(G)$, this defines an Euler characteristic in $R(G)$. The Euler characteristic in $K$-theory has an alternate description $\chi^{K_G}(X)=\sum (-1)^i[H^i(M,\C)]$ by considering the cohomology groups $H^i(M,\C)$ as classes in $R(G)$ \cite[Remark V.1.15]{LMS86}. 
If a statement is true about both $\chi^G(M)$ and $\chi^{K_G}(M)$ without leveraging specific properties of either $\Sphere$ or $KU_G$, we write $\chi^G(M)$.

We say $8$ points of $\CP^2$ are in general position if they belong to the open subset of $(\CP^2)^8$ consisting of points such that the set of cubic polynomials passing through them has dimension $2$ and the corresponding pencil of cubics passing through them has fibers with at worst nodal singularities. 

\begin{theorem}\label{thm:rational_cubic_count}
 Let $G$ be a finite group acting on $\CP^2$, and let $S$ be a $G$-invariant set of $8$ points in general position in $\CP^2$. 
    There exists a formula $N_{3, \CP^2, S}^G$ in $A(G)$ and $R(G)$ that counts orbits of rational plane cubics through $S$ in $\CP^2$, and \begin{equation}\label{eq:rational_cubic_count}N_{3,\CP^2,S}^G=\sum_{\substack{ G\cdot C,\\ C \text{ a rational cubic }\\ \text{passing through }S}} \tr_{G_{C}}^G m^{G_{C}}(C), 
    \end{equation}
    where $G_C$ is the stabilizer of the rational cubic $C$. 
\end{theorem}

First,  observe the set of rational cubics through $S$ is a finite $G$-set. 
 \begin{lemma}\label{lem:count_is_invariant}
    The set of 12 rational cubics through $S$ is a $G$-invariant subset of $\CP^2$. 
\end{lemma}
\begin{proof}
    Note that cubic curves are either smooth, elliptic curves or rational and singular. 
    Since the points of $S$ are general, the rational cubics through $S$ are nodal with 1 node, so we need to show that $G\cdot C$ contains only nodal cubics passing through $S$ for any rational cubic $C = \text{im}(u)$, $u\colon \CP^1\to \CP^2$ a degree $3$ rational curve. The group action $G\times \CP^2\to \CP^2$ is a function so must act transitively on the set of nodal points of rational cubics within the same orbit, and in particular no rational, nodal cubics share an orbit with smooth cubics. All rational cubics in an orbit must pass through $S$ since $G\cdot S=S$. 
    \end{proof}

Thus there is a class representing $N_{3, \CP^2,S}^G$ in $A(G)$ as a $G$-set and a class in $R(G)$ under the map $\pi_0^G\Sphere\to \pi_0^GKU_G$. To derive the explicit formula in Theorem \ref{thm:rational_cubic_count}, we translate the question of counting $G$-orbits of rational cubics through $S$ to a question of counting nodal cubics in a $G$-invariant general pencil of cubics in $\CP^2$. See \cite[Section 9.2]{KLSW-rational-curves} for a non-equivariant version of this construction. 

Let \[S:=\sum_{(H_i)\subset G} n_i[G/H_i],\] with $n_i\in\Z_{\geq 0}$, be representative class in $A(G)$ of a $G$-invariant set of $8$ general points in $\CP^2$. A choice of basis 
for $\op{H}^0(\CP^2, -K_{\CP^2})$ determines an embedding $\CP^2\hookrightarrow \CP^9$. 
There are 8 linear conditions on $\op{H}^0(\CP^9, \mathcal{O}(1))$ corresponding to vanishing at the 8 points in $S$, and these conditions are linearly independent since the points in $S$ are in general position. This implies  that     
\begin{equation}\label{eq:big_space}
\{F\in \op{H}^0(\CP^9,\mathcal{O}(1))\colon F(p_i)=0\text{, }\forall p_i\in S\}
\end{equation} 
is 2-dimensional. Let $\langle f, g\rangle$ be a basis of this vector space. Note that $\langle f, g\rangle$ is a $G$-invariant subspace of \eqref{eq:big_space}. 

Let $\{sf(p)+tg(p)=0\colon [s,t]\in \CP^1\}\subseteq \CP^2$ be the pencil spanned by $f$ and $g$, and let 
\[
X:=\{([s,t],p)\colon sf(p)+tg(p)=0\}\subseteq \CP^1\times \CP^2
\] be the total space of the pencil. Since $\langle f, g\rangle$ is $G$-invariant, $X$ is $G$-invariant by construction, as orbits of cubics through $S$ must only contain other cubics passing through $S$. 
Also note that $X$ is a general pencil as $S$ is a general set of $8$ points, whence $f$ and $g$ pass through a $9^{\text{th}}$ point $P$ by Cayley-Bacharach such that $S\cup \{P\}$ is a set of 9 general points. Let $B:=S\cup \{P\}$ be the base locus of $X$. By construction, $B$ is a genuine, finite $G$-set. Furthermore, 
\[B=S+\{*\}=\sum_{(H_i)\subset G} n_i[G/H_i]+\{*\},
\]$n_i\in\Z_{\geq 0}$, is a representative for the class of $B$ in $A(G)$ in terms of transitive $G$-sets. 

There are two projections from $X$,
\begin{equation}
    \begin{tikzcd}\label{diag:pencil_projections}
        & X\arrow["\pi_2"]{dr}\arrow[swap, "\pi_1"]{dl} & \\
        \CP^1 & & \CP^2
    \end{tikzcd}
\end{equation}
For any $p\in \CP^1$, we will write $X_p$ for the cubic in $X$ which is the fiber of $\pi_1$ over $p$. The orbits of fibers of $\pi_1$ over the set \[
\Sing(\pi_1):=\{p\colon \pi_1^{-1}(p) \text{ is nodal}\} \subseteq \CP^1
\] describe the orbits of rational cubics through $S$. More specifically, any $X_p=\pi_1^{-1}(p)$ is either a smooth, elliptic cubic if $p\in \CP^1\setminus\Sing(\pi_1)$ or a rational, nodal cubic if $p\in \Sing(\pi_1)$. Note too that nodal fibers $X_p$ are exactly the zeros of $d\pi_1$. 
We deduce
\[
\bigl\{\text{orbits of rational cubics through }S  \bigr\} \stackrel{\sim}{\longleftrightarrow}\bigl\{\text{orbits of nodal cubics in }X  \bigr\}. 
\]
In particular, $N_{3,\CP^2, S}^G$ can be computed in terms of orbits of nodal cubics over $\Sing(\pi_1)$ in $X$.  
We will ultimately show the count of orbits of nodal cubics in $X$ is equal to $N_{3, \CP^2, S}^G$ in $A(G)$ and $R(G)$ (Corollary \ref{cor:final_count}). 

\begin{theorem}\label{thm:pencil_count}
    Let $G$ be a finite group acting on $\CP^2$, and let $X$ be a $G$-invariant general pencil of cubic curves in $\CP^2$ with base locus $B=\sum_{(H_i)\subset G} n_i[G/H_i]+\{*\}$ in $A(G)$. Then the count of orbits of nodal cubics in $X$ is given by the Euler characteristic $\chi^G(X)$,
    \begin{equation}\label{eq:count_of_cubics}
    \chi^G(X) = \sum_{\substack{ G\cdot X_p,\\ X_p\in X \text{ nodal}}} \tr_{G_{X_p}}^G m^{G_{X_p}}(X_p)
    \end{equation}
    in $A(G)$ and $R(G)$. When $G$ is additionally an abelian group and we write $\CP^2\cong \mathbb{P}V$ for some 3-dimensional $G$-representation $V$, $\chi^{K_G}(X)$ can be computed as: 
\begin{equation}\label{eq:projective_bundle_formula_for_pencil}{\bigwedge}^{\raisebox{-0.4ex}{\scriptsize$2$}} \text{ }V+\sum_{\substack{G\cdot p_i,\\p_i\in B}}\tr_{G_{p_i}}^G(-1+V|_{H_i})=\sum_{\substack{ G\cdot X_p,\\ X_p\in X \text{ nodal}}} \tr_{G_{X_p}}^G m^{K_{G_{X_p}}}(X_p)
    \end{equation} in $R(G)$, where $H_i=\op{stab}(p_i)$. 
\end{theorem}

This statement is quite general, as we do not impose any hypotheses on the group action (e.g., free, algebraic), simply that the group is finite and, for \eqref{eq:projective_bundle_formula_for_pencil} specifically, abelian.  

We first derive a projective bundle formula needed for the second part of Theorem \ref{thm:pencil_count} using formal group laws for complex equivariant $K$-theory in the case of an abelian group. Note that we implicitly use Theorem \ref{thm:gauss_bonet}.

\begin{lemma}\label{lem:proj_bundle_formula}
Let $G$ be a finite, abelian group and let $V$ be an $(n+1)-$dimensional $G$-representation. Then 
\begin{equation}\label{eq:proj_bundle_equation}
    \chi^{K_G}(\mathbb{P}V)={\bigwedge}^{\raisebox{-0.4ex}{\scriptsize$n$}} \text{ }V
\end{equation} in $R(G)$. 
\end{lemma} 
\begin{proof}
    We have a short exact sequence on $\mathbb{P}V$,
    \[ 0\to \mathcal{O}_{\mathbb{P}V}(-1)\to V\otimes \mathcal{O}_{\mathbb{P}V}\to Q\to 0
    \]
    where $Q$ is the quotient bundle. Tensoring with $\mathcal{O}_{\mathbb{P}V}(1)$ gives 
    \begin{equation}\label{seq:for_chern}
    0\to\mathcal{O}_{\mathbb{P}V} \to V\otimes \mathcal{O}_{\mathbb{P}V}(1)\to T\mathbb{P}V\to 0. 
    \end{equation}
    Let $c$ denote the total chern class in equivariant $K$-theory. Exactness of \eqref{seq:for_chern} implies 
    \begin{align*}
        c(V\otimes \mathcal{O}_{\mathbb{P}V}(1)) &= c(\mathcal{O}_{\mathbb{P}V})c(T\mathbb{P}V)\\
        &=c(T\mathbb{P}V). 
    \end{align*}

    We can write $V\cong \oplus_{i=0}^n L_i$ as a sum of line bundles $L_i$ corresponding to 1-dimensional $G$-representations by the splitting principle in equivariant $K$-theory. 
    Then
    \begin{align*}
        c(T\mathbb{P}V)&= c(V\otimes \mathcal{O}_{\mathbb{P}V}(1))\\
        &=c(\oplus_{i=0}^n(L_i\otimes \mathcal{O}_{\mathbb{P}V}(1)) \\ 
        &= \Pi_{i=0}^n(1+c_1(L_i\otimes\mathcal{O}_{\mathbb{P}V}(1))) \\
        &= \Pi_{i=0}^n (1+c_1(L_i)+c_1(\mathcal{O}_{\mathbb{P}V}(1))-\beta c_1(L_i)c_1(\mathcal{O}_{\mathbb{P}V}(1))
    \end{align*}
    where $\beta$ is the Bott element for the formal group law in equivariant complex $K$-theory. The study of equivariant formal group laws is rich, see \cite{CGK-fgl}, \cite{Greenlees-fgl}, and \cite{Hausmann-fgl}.

    Let $\alpha$ denote $c_1(\mathcal{O}_{\mathbb{P}V}(1))$ to simplify notation. Recall using the Atiyah-Hirzebruch spectral sequence that
    \[
    K_G^*(\mathbb{P}V)\cong K_*^G\oplus \alpha K_*^G \oplus \cdots\oplus \alpha^n K_*^G. 
    \]
    The pushforward to a point $\mathbb{P}V\stackrel{p}{\to} pt$ induces a $K_*^G$-module homomorphism $p_*$ which reads off the coefficient of the $\alpha^n$ term. Rearranging terms of $c(T\mathbb{P}V)$, the $\alpha^n$ coefficient of $c(T\mathbb{P}V)$ is
    \begin{equation}\label{eq:random1}
    \left(\sum_{i=0}^n(1+c_1(L_i)) \right)\prod_{j\ne i} (1-\beta c_1(L_i)). 
    \end{equation}
    Since $c_1(L_i)=\beta^{-1}(1-L_i)$, we can write  $\Pi_{j\ne i} (1-\beta c_1(L_i)) 
    = \Pi_{j\ne i} L_i$. Thus \eqref{eq:random1} simplifies to 
    \begin{align*}
        \left(\sum_{i=0}^n(1+c_1(L_i)) \right)\prod_{j\ne i} (1-\beta c_1(L_i)) &= \sum_{i=0}^n \prod_{j\ne i} L_i\\
        &= {\bigwedge}^{\raisebox{-0.4ex}{\scriptsize$n$}} \text{ }V. 
    \end{align*}
    We conclude 
    \[
    \chi^{K_G}(\mathbb{P}V) = p_* (c(T\mathbb{P}V)) = {\bigwedge}^{\raisebox{-0.4ex}{\scriptsize$n$}} \text{ }V, 
    \]
    as desired. 
\end{proof}

\noindent\textit{Proof of Theorem \ref{thm:pencil_count}.} We will compute $\chi^G(X)$ in two ways using both of the projections in \eqref{diag:pencil_projections} and set them equal to prove Theorem \ref{thm:pencil_count}. 
First we consider $\pi_1\colon X\to \CP^1$. As before let 
\[
\Sing(\pi_1):=\{p\colon \pi_1^{-1}(p) \text{ is nodal}\}\subseteq \CP^1
\]
so that $\CP^1=(\CP^1-\Sing(\pi_1))\amalg \Sing(\pi_1)$. 
Note that $\chi^G$ is multiplicative in products, hence in fiber bundles. This implies that  $$\chi^G(X) = \chi^G(\CP^1\setminus\Sing(\pi_1))\chi^G(X_\eta)+\chi^G(\Sing(\pi_1))\chi^G(X_p),$$ where $X_\eta$ is a general smooth fiber and $X_p$ is a rational, nodal fiber over $p\in \Sing(\pi_1)$. 
Note also that $\chi^{top}(X_\eta)=0$ non-equivariantly, whence $\chi^{G_{X_\eta}}(X_{\eta})=0$. We deduce $$\chi^G(\CP^1\setminus\Sing(\pi_1))\cdot\chi^{G_{X_\eta}}(X_{\eta})=0.$$ Thus
\begin{align*}
    \chi^G(X) &= \chi^G(\Sing(\pi_1))\chi^G(X_p)\\
    &=\sum_{\substack{G\cdot X_p,\\ X_p\in X\text{ nodal}}} 
    \tr_{G_{X_p}}^G\chi^{G_{X_p}}(X_p).\\
\end{align*}

This technically completes the first part proof. We will compute $\chi^G(X)$ again using $\pi_2\colon X\to \CP^2$. Note that $X\cong \op{Bl}_B\CP^2$. The fact that this isomorphism exists is a classical fact, and this isomorphism is equivariant because $B$ is a $G$-invariant subset of $\CP^2$. 
Thus
\begin{align*}
\chi^G(X)&=
\chi^G(\CP^2)-\chi^G(B)+\chi^G(\mathbb{P}N_B\CP^2),
\end{align*}
noting that $\mathbb{P}N_B\CP^2\subset X$ is the exceptional divisor.

We deduce that 
\begin{equation}\label{eq:nodes_equals_euler}
\chi^G(\CP^2)-\chi^G(B)+\chi^G(\mathbb{P}N_B\CP^2) = \sum_{\substack{G\cdot X_p,\\ X_p\in X\text{ nodal}}}\tr_{G_{X_p}}^G m^{G_{X_p}}(X_p)
\end{equation}
in both $A(G)$ and $R(G)$. 

Equation \eqref{eq:nodes_equals_euler} holds when $G$ is any finite group. When $G$ is both finite and abelian, we use Lemma \ref{lem:proj_bundle_formula} to prove the second part of Theorem \ref{thm:pencil_count}. Writing $\CP^2\cong \mathbb{P}V$ for some $G$-representation $V$, the left hand side of equation \eqref{eq:nodes_equals_euler} of $\chi^{K_G}(Bl_B\mathbb{P}V)$ is
\[
\chi^{K_G}(\mathbb{P}V)-\chi^{K_G}(B)+\chi^{K_G}(\mathbb{P} N_B\mathbb{P}V) 
\]
in $R(G)$. 

Let $T_{p_i}\mathbb{P}V$ denote the tangent space of $\mathbb{P}V$ at $p_i\in B$, which is a $G_{p_i}$-representation. Recall that 
\[B=\sum_{(H_i)\subset G} n_i[G/H_i]+\{*\},
\]
so each $p_i$, which is a representative of its orbit in $B$, has stabilizer $H_i$. Then 
\[
\chi^{K_G}(\mathbb{P} N_B\mathbb{P}V)=\sum_{\substack{G\cdot p_i,\\ p_i\in B}}\tr_{H_i}^GT_{p_i}V, 
\]
whence 
\begin{align*}
\chi^{K_G}(Bl_B\mathbb{P}V) &= \chi^{K_G}(\mathbb{P}V)-\chi^{K_G}(B)+\chi^{K_G}(\mathbb{P} N_B\mathbb{P}V)\\
&= {\bigwedge}^{\raisebox{-0.4ex}{\scriptsize$2$}} \text{ }V+\sum_{\substack{G\cdot p_i\\ p_i\in B}}\tr_{H_i}^G(-1+T_{p_i}\mathbb{P}V)
\end{align*}
in $R(G)$ using Lemma \ref{lem:proj_bundle_formula} and the fact that $\chi^G(\mathbb{P}N_B(\mathbb{P}V))$ is a disjoint union of $T_{p_i}\mathbb{P}V$ over the orbit type of $B$. Additionally,
\[
N_B \mathbb{P}V \cong \mathcal{O}(3)\vert_B \oplus \mathcal{O}(3)\vert_B.
\]
Since $B$ is a disjoint union of points, $\mathcal{O}(3) \cong \mathcal{O}$ locally. Further, generically we have that all $p_i$ are in $V \subset \mathbb{P}(V)$. Thus $T_{p_i} \mathbb{P}(V) \cong T_{p_i} V \cong V\vert_{H_i}$. 
 We conclude that equation \eqref{eq:nodes_equals_euler} is 
\begin{equation}\label{eq:nodes_equals_euler_projective}{\bigwedge}^{\raisebox{-0.4ex}{\scriptsize$2$}} \text{ }V+\sum_{\substack{G\cdot p_i\\ p_i\in B}}\tr_{H_i}^G(-1+V|_{H_i}) = \sum_{\substack{G\cdot X_p,\\ X_p\in X\text{ nodal}}}\tr_{G_{X_p}}^G m^{G_{X_p}}(X_p) \end{equation}
in $R(G)$ when $G$ is finite and abelian, as desired. 
\hfill $\Box$

Singling out equations \eqref{eq:nodes_equals_euler} and \eqref{eq:nodes_equals_euler_projective} in the proof, we have shown:

\begin{corollary}\label{cor:final_count} If $G$ is a finite group acting on $\CP^2$ and $S$ is a $G$-invariant set of 8 general points in $\CP^2$, an equivariantly enriched count of rational cubics through $S$ is given by
    \begin{align*}
    N_{3, \CP^2, S}^G &= \sum_{\substack{ G\cdot C,\\ C \text{ a rational cubic }\\ \text{passing through }S}} \tr_{G_p}^G m^{G_C}(C) \\
    &= \chi^G(\CP^2)-\chi^G(B)+\chi^G(\mathbb{P}N_B\CP^2) 
    \end{align*}
in $A(G)$ and $R(G)$, where $B$ is the base locus of the $G$-invariant pencil of cubics the rational cubics through $S$ interpolate. When $G$ is both finite and abelian, 
    \begin{align*}
    N_{3, \CP^2, S}^G &= \sum_{\substack{ G\cdot C,\\ C \text{ a rational cubic }\\ \text{passing through }S}} \tr_{G_p}^G m^{K_{G_C}}(C)  \\
    &= {\bigwedge}^{\raisebox{-0.4ex}{\scriptsize$2$}} \text{ }V+\sum_{\substack{G\cdot p_i\\ p_i\in B}}\tr_{H_i}^G(-1+V|_{H_i})
    \end{align*}
    in $R(G)$, where $H_i$ is the stabilizer of $p_i$. The cardinality (resp. rank) of $N_{3, \CP^2, S}^G$ in $A(G)$ (resp. in $R(G)$) is 12, recovering the non-equivariant count of 12 rational cubics through $S.$
\end{corollary}

By taking $H$-fixed points of $N_{3, \CP^2, S}^G\in A(G)$ for any $H\leq G$, we obtain an integer count of the number of orbits of rational cubics stabilized by $H$. Similarly, we may look at group characters of $N_{3,\CP^2, S}^G\in R(G)$. In general, the map $\pi_0^G \Sphere\to \pi_0^G KU_G$ induced from the unit $\Sphere\to KU_G$ is not injective, so we cannot always expect to obtain a unique formula for $N_{3, \CP^2, S}^G$ in $A(G)$ from one in $R(G)$.  

As mentioned at the beginning of this section, Theorem \ref{thm:rational_cubic_count} is also related to the question of counting real rational cubics through a size $8=n_1+2n_2$ set $S$ of general points in $\CP^2$ with $n_2$ is the number of complex conjugate pairs of points in $S$ and $n_1$ is the number of real points. If $G=\Z/2$ acts on $\CP^2$ by pointwise complex conjugation, then \[S=n_1[G/G]+n_2[G]\] is the orbit type of $S$ in $A(G)$. The character of the 12-dimensional representation $N_{3,\CP^2, S}^G$ considered in $R(G)$ will recover the signed count of real rational cubics through $S$, likewise the number of $G$-fixed points of $N_{3, \CP^2, S}^G$ considered in $A(G)$ will recover the signed count. 

The question of counting real rational curves was studied by Welschinger in \cite{welschinger-1}, \cite{welschinger-2} by assigning a sign to each real nodal point of a rational curve  to obtain a signed count of real rational curves. Given any real rational curve $u\colon S^2\to \CP^2$ with nodal points $p_1, \ldots, p_n$, Welschinger defines the \textit{mass of $u$ at $p_i$, $m_{p_i}(u)$}, to be  $+1$ if $p_i\in\mathbb{RP}^2$ defines a split node of $\op{im}(u)$, $-1$ if $p_i\in\mathbb{RP}^2$ defines a non-split node of $\op{im}(u)$, and any $p_i\in \mathbb{CP}^2\setminus \mathbb{RP}^2$ which occur in a complex conjugate pair of nodal points are not counted. By definition, $u$ is \textit{split} at $p_i\in \mathbb{RP}^2$ if $\text{im}(u)$ has two real branches at $p_i$, and $u$ is \textit{non-split} at $p_i\in \mathbb{RP}^2$ if $\op{im}(u)$ has a complex conjugate pair of branches at $p_i$. 
The Welschinger mass of $u$ is then defined to be 
\begin{equation}\label{eq:welch_mass}
m(u):= (-1)^{\#\text{ non-split nodes}} = \prod_{\substack{p \text{ a real} \\ \text{node of } \op{im}(u)}} {m_p(u)}.
\end{equation}
Since rational cubics (generically) have exactly one ordinary double point $p$, the Welschinger mass of a real rational cubic is simply $m(u)=+1$ if $u$ is split at $p$ and $m(C)=-1$ if $C$ is non-split at $p$. 

\begin{figure}[ht]\label{figure:node_types}
\centerline{\includegraphics[width=.5\textwidth]{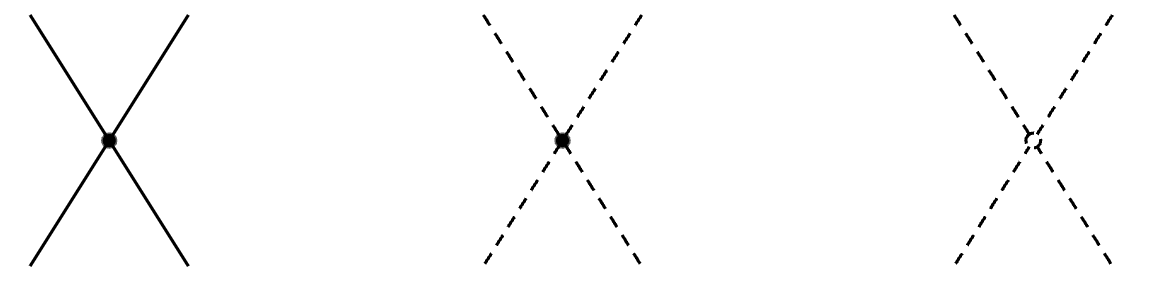}}
\caption{Real split, real non-split, and complex conjugate node branches.}
\end{figure}

A significant accomplishment of Welschinger in real Gromov-Witten theory is given by the signed count of real rational curves:

\begin{theorem}\label{thm:original_welschinger}
(Welschinger \cite{welschinger-1, welschinger-2}) Fix integers $d, n_1$, and $n_2$ such that $n_1+2n_2=3d-1$. For any generic choice of $n_1$ real and $n_2$ complex conjugate pairs of points in $\CP^2$, the sum 
\[
W_{d, n_1}=\sum_{\substack{u\text{ degree d } \\\text{real rational curve} \\\text{through the points}}} m(u)
\]
is independent of choice of $3d-1=n_1+2n_2$ points. 
\end{theorem}

 Work of J. Solomon computes $W_{d, n_2}$ recursively, see \cite{solomon-thesis}. 
 Recent work of J. Kass, M. Levine, J. Solomon, and the second author of this work introduces quadratically enriched Gromov-Witten invariants,  effectively defining a quadratically enriched count of rational curves through a suitable divisor on a del Pezzo surface $X$ of degree $\geq 4$ over fields of characteristic $\neq 2, 3$ \cite{KLSW-rational-curves}. This builds on prior work of M. Levine on quadratically enriched Welschinger invariants \cite{levine-welschinger}. 

In a particular case of Welschinger's work, the signed count of real rational plane cubics through  $8=n_1+2n_2$ general points in $\CP^2$ is $$W_{3, n_1}=n_1=8-2n_2.$$ This is discussed in  \cite[Example 1]{welschinger-2}. We show this equivariantly in Corollary \ref{cor:real_count}.  

\begin{remark} 
It's worth noting that the genuine, unsigned count of real rational cubics through 8 points (non-equivariantly) is not the same as the signed count given by $W_{3, n_2}$. For example, when $S$ contains 8 real general points, then the signed count of real rational cubics through $S$ is $W_{3,0} = 8$. However, the genuine number of real rational cubics through $S$ is between $8$ and $12$, and there are choices of the 8 points for which both $8$ and $12$ real rational cubics occur, see \cite[Proposition 4.7.3]{DK-2000}. See \cite[7.2.2]{Sottile-real} for an illuminating discussion of the question of how many real rational cubics interpolate any 8 of the 9 real points in a pencil and an example of a pencil in which all 12 rational cubics are real. Separately, see also \cite[3.1]{IKS-2004} for an example of a pencil where 8 of the 9 base locus points are complex.
\end{remark}

The equivariant mass in Definition \ref{defn:mass} is an equivariant enrichment of the mass used by Welschinger in \cite{welschinger-1, welschinger-2} and the equivariant weight of a node defined by the first author in \cite{betheapencil}. If $G=\Z/2$ acts on $\CP^2$ by pointwise complex conjugation and $X_p=\pi_1^{-1}(p)$, $p \in \Sing(\pi_1)$, is a nodal, rational cubic through the $G$-invariant set $$S=n_1[G/G]+n_2[G]$$ with stabilizer $G_p$, then the transfer of the mass of $X_p$, $m^{G_p}(X_p)=\chi^{G_p}(X_p)$, is 
\begin{equation}\label{eq:mass_in_A(G)}
\tr_{G_p}^G m^{G_p}(X_p)=
\begin{cases} 
      [G]-[G/G], & \text{if $X_p$ is real, split} \\
      [G/G] & \text{if $X_p$ is real, non-split} \\
      [G] & \text{if the nodal point of $X_p$ is complex}\\
   \end{cases}
\end{equation}
in $A(G)$ and 
\begin{equation}\label{eq:mass_in_R(G)}
\tr_{G_p}^G m^{G_p}(X_p)=
\begin{cases} 
      1_{\sigma}, & \text{if $X_p$ is real, split} \\
      1& \text{if $X_p$ is real, non-split} \\
      1+1_{\sigma} & \text{if the nodal point of $X_p$ is complex}\\
   \end{cases}
\end{equation}
in $R(G)$, where $1$ denotes the trivial representation and $1_{\sigma}$ is the sign representation. This can be seen directly from a $G$-CW presentation of each rational curve type. Note when the nodal point of $X_p$ is real, whether split or non-split, $G_p=G$ and $\tr_{G_p}^G m^{G_p}(X_p) = m^G(X_p)$. When the nodal point of $X_p$ is complex, 
$G_p=e$ and $\tr_{G_p}^G m^{G_p}(X_p) = \tr_e^G [G/G] = [G]$. 
    
\begin{corollary}\label{cor:real_count} 
     Let $G=\Z/2$ act on $\CP^2$ by complex conjugation, and let $S=n_1[G/G]+n_2[G]$ be a $G$-invariant set of $8$ points in $\CP^2$. Then 
    \[
    -|(N_{3, \CP^2, S}^G)^G|=W_{3, n_1}
    \]
    when $N_{3, \CP^2, S}^G$ is considered in $A(G)$, and
    \[
    -\chi_{N_{3, \CP^2, S}^G}(g) 
    = W_{3, n_1}
    \]
    when $N_{3, \CP^2, S}^G$ is considered in $R(G)$. 
    Thus $-N_{3, \CP^2, S}^G$ recovers the $\Z$-valued signed count of real rational cubics through $S$, $W_{3, n_1}$, when $G=\Z/2$ acts on $\CP^2$ by pointwise complex conjugation. 
\end{corollary}
The notation $\chi_{N_{3, \CP^2, S}^G}(g)$  
denotes the character of the representation $N_{3, \CP^2, S}^G$ 
for the non-identity element $g$ of $G$. 
\begin{proof}
Taking fixed points of \eqref{eq:mass_in_A(G)}, we see
\begin{equation}\label{eq:fixed_points_of_A(G)_mass}
|\tr_{G_p}^G m^{G_p}(X_p)|=
\begin{cases} 
      -1, & \text{if $X_p$ is real, split} \\
      1 & \text{if $X_p$ is real, non-split} \\
      0 & \text{\text{if the nodal point of $X_p$ is complex}}\\
   \end{cases}
\end{equation}
in $\Z$. Likewise, the character for the non-identity element of $G$ of \eqref{eq:mass_in_R(G)} gives 
\begin{equation}\label{eq:fixed_points_of_R(G)_mass}
\chi_{\tr_{G_p}^G m^{G_p}(X_p)}(g)=
\begin{cases} 
      -1, & \text{if $X_p$ is real, split} \\
      1 & \text{if $X_p$ is real, non-split} \\
      0 & \text{\text{if the nodal point of $X_p$ is complex}}\\
   \end{cases}
\end{equation}
in $\Z$. In particular, we recover the Welschinger mass from both 
\[
-|\tr_{G_p}^G m^{G_p}(X_p)| = m(X_p)
\]
and 
\[
-\chi_{\tr_{G_p}^G m^{G_p}(X_p)}(g) = m(X_p). 
\]
We deduce immediately that 
\begin{align*}
-|(N_{3, \CP^2, S}^G)^G| &= \sum_{\substack{ G\cdot X_p,\\ X_p\in X \text{ nodal}}} -|\tr_{G_{p}}^G m^{G_{X_p}}(X_p)| \\
&= \sum_{\substack{X_p\text{ rational } \\ \text{through $S$}}} m(X_p)\\
&=  W_{3, n_1}. 
\end{align*}

When we consider $N_{3, \CP^2, S}^G$ as an element of the representation ring, we likewise conclude
\begin{align*}
-\chi_{N_{3, \CP^2, S}^G}(g)&=\sum_{\substack{ G\cdot X_p,\\ X_p\in X \text{ nodal}}}  -\chi_{\tr_{G_{X_p}}^G m^{G_{X_p}}(X_p)}(g)\\
&= \sum_{\substack{X_p\text{ rational } \\ \text{through $S$}}} m(X_p)\\
&=  W_{3, n_1}. 
\end{align*}
\end{proof}

We've shown that 
\begin{equation}\label{eq:wel_euler_char_1}
|\chi^G(X^G)|=\chi_{\chi^G(X^G)}(g)=-W_{3, n_1}=-n_1, 
\end{equation}
where $X$ is still the $G$-invariant pencil with base locus 
\[
[B]=S +\{*\} = (n_1+1)[G/G]+n_2[G]
\]
in $A(G)$ whose orbits of nodal cubics are in bijection with orbits of rational cubics through $S$. The blowup formula for $\chi^G(X)$ also shows  \eqref{eq:wel_euler_char_1} holds. Observe
\begin{align*}
    \chi^G(X) & = \chi^G(\CP^2)+\chi^G(B)\chi^G(\CP^1)-\chi^G(B) \\
    &= \chi^G(\CP^2) +(n_1+1)[G/G]\cdot\chi^G(\CP^1)+n_2[G]\cdot\chi^G(\CP^1) - \chi^G(B). 
\end{align*}
Taking fixed point cardinalities of $\chi^G(X^G)$ gives
\begin{align*}
    |\chi^G(X^G)| &= |\chi^G(\mathbb{RP}^2)+(n_1+1)\chi^G(\mathbb{RP}^1)-\chi^G(B^G)| \\ 
    &= -n_1, 
\end{align*}
another way to see 
\[
-|(N_{3, \CP^2, S}^G)^G| = -|\chi^G(X)| = n_1 = W_{3,n_1}. 
\]
It's clear by the same reasoning that $-\chi_{\chi^G(X)}(g)=W_{3, n_1}$. Furthermore, $N_{3, \CP^2, S}^G = \chi^G(X)$ has a projective bundle formula in $R(G)$, which is 
\[
N_{3, \CP^2, S}^G = {\bigwedge}^{\raisebox{-0.4ex}{\scriptsize$2$}} \text{ }V+\sum_{\substack{G\cdot p_i,\\p_i\in B}}\tr_{G_{p_i}}^G(-1+V|_{H_i})
\]
where we now write $V=\mathbb{C}^3$ with complex conjugation action so that $\mathbb{P}V=\CP^2$ inherits a coordinate-wise complex conjugation action and $H_i$ is the stabilizer of $p_i$. Computing $-\chi_{\chi^G(X)}(g)=W_{3, n_1}$ with the projective bundle formula also recovers $W_{3, n_1}$. 

We've shown $N_{3, \CP^2, S}^G$ encodes orbit information about the 12 rational cubics through $S$, recovers the count of 12 rational cubics through $S$, and recovers the Welschinger invariant by taking conjugation fixed points. In other words, the equivariant topology of the invariant pencil the 12 rational cubics interpolate determines both the complex and real signed count of rational cubics in a single equivariant formula. 

\begin{remark} 
    Formula \eqref{eq:mass_in_A(G)} for $\tr_{G_p}^G m^{G_p}(X_p)$ in the $\Z/2$ case was computed directly by writing down a $G$-CW complex for each type of rational cubic in $X$. 
    Computing just one of \eqref{eq:mass_in_A(G)} or \eqref{eq:mass_in_R(G)} is enough to obtain both, as $A(G)\cong R(G)$ since $G$ is a finite cyclic group \cite[Proposition 4.5.4]{TD79}.
\end{remark}

\section{An Explicit Euler Number Computation}\label{section:Hom_bundle}

Section \ref{section:example} shows 
\[
N_{3, \CP^2, S}^G = \chi^G(X)
\]
where $X$ is the general pencil of plane cubics with base locus $B$ represented by $[S]+\{*\}$ in $A(G)$. By Theorem \ref{thm:gauss_bonet}, $\chi^{K_G}(X) = {\pi_X}_* e^{K_G}(TX)$ where $\pi_X\colon X\to *$ is the map to a point in the following diagram
\begin{equation}
    \begin{tikzcd}
            X\cong \op{Bl}_B\bbP V_2 \arrow["\pi_2"]{r}\arrow[swap, "\pi_1"]{d}\arrow["\pi_X"]{dr} & \bbP V_2\supset B \arrow["p_2"]{d} \\ 
        \bbP V_1 \arrow[swap,"p_1"]{r} & * \\
    \end{tikzcd}
\end{equation}
and we've written $\CP^2 = \mathbb{P}V_2$ and $\CP^1 = \mathbb{P}V_1$ for rank 3 and 2 complex $G$-representations $V_2$ and $V_1$ respectively. In this section we separately compute the Euler number of a Hom-bundle with respect to a section, which also computes the number of nodal cubics in $X$. Let $W:=\op{Hom}(\pi_1^*T^*\bbP V_1, T^*X) \to X$, which is a rank 2 bundle over $X$. Note $d\pi_1\colon X\to W$ is a section. 
We compute $\pi_{X*}e^{K_G}(W) = n_G(W,d\pi_1)$, and thus the $R(G)$-enriched equivariant count of nodal cubics in $X$, explicitly. 

We have the short exact sequence on X
\begin{equation}\label{hom_eq1}
0\to \pi_1^*\mathcal{O}(-1)\otimes \pi_2^*\mathcal{O}(-3)\to \pi_1^*T^*\bbP V_1|_X\oplus \pi_2^*T^*\bbP V_2|_X\to T^*X\to 0. 
\end{equation}
Note that $\pi_{X*}e^{K_G}(W)$ is represented by applying $\pi_{X*}$ to 
\[
0 \to \wedge^2W^* \to W^* \to \mathcal{O}_X\to 0,
\]
so that 
\begin{equation}\label{hom_eq2}
    \pi_{X*} e^{K_G}(W) = R\pi_{X*}\mathcal{O}_X + R\pi_{X*}\det W^*-R\pi_{X*}W^*.
\end{equation}
Applying $\op{Hom}(\pi_1^*T^*\bbP V_1,-)$ to \eqref{hom_eq1} and dualizing, we have 
\begin{equation}\label{hom_eq2.5_dualized}
    0\to W^*\to \mathcal{O}_X\oplus(\pi_1^*\mathcal{O}(-2)\otimes \pi_2^*T\bbP V_2)\to \pi_1^*\mathcal{O}(-1)\otimes \pi_2^*\mathcal{O}(3)\to 0, 
\end{equation}
thus 
\[
R\pi_{X*}W^* = R\pi_{X*} \mathcal{O}_X + R\pi_{X*}(\pi_1^*\mathcal{O}(-2)\otimes \pi_2^* T\bbP V_2) - R\pi_{X*}(\pi_1^*\mathcal{O}(-1)\otimes \pi_2^*\mathcal{O}(3). 
\]
Combining this with \eqref{hom_eq2},  
\begin{equation}\label{hom_eq4}
    \pi_{X*}e^{K_G}(W) = R\pi_{X*}\det W^* + R\pi_{X*}(\pi_1^*\bbO(-1)\otimes \pi_2^* \bbO(3)) - R\pi_{X*}(\pi_1^*\bbO(-2)\otimes \pi_2^* T\bbP V_2). 
\end{equation}

Note that \eqref{hom_eq2.5_dualized} also implies 
\begin{equation}\label{hom_eq5}
\det W^* \cong \pi_1^*\bbO(-3)\otimes \pi_2^*\bbO(-3)\otimes \det \pi_2^*T\bbP V_2. 
\end{equation}
From the Euler sequence on $\bbP V_2$, 
\[
0\to \bbO_{\bbP V_2} \to \bbO_{\bbP V_2}(1)\otimes V_2\to T\bbP V_2\to 0,
\]
we see 
\begin{equation}\label{hom_eq8}
\det T\bbP V_2 \cong \bbO_{\bbP V_2}(3)\otimes \det V_2.\end{equation}
Combining \eqref{hom_eq5} and \eqref{hom_eq8},
\begin{equation}\label{hom_eq9}
    \det W^* \cong \pi_1^*\bbO(-3)\otimes \pi_2^*\det V_2.
\end{equation}
From the Euler sequence, we also see $\pi_2^*T\bbP V_2 \cong \pi_2^*(\bbO(1)\otimes V_2)-\pi_2^*\bbO(-2)$. 

Putting this together, we have 
\begin{align}\label{eq_with_4_terms}
    \pi_{X*}e^{K_G}(W)&= R\pi_{X*} \det W^* + R\pi_{X*}(\pi_1^*\bbO(-1)\otimes \pi_2^*\bbO(3)) - R\pi_{X*}(\pi_1^*\bbO(-2)\otimes \pi_2^*T\bbP V_2) \notag \\
    &= R\pi_{X*}(\pi_1^*\bbO(-3)\otimes \pi_2^*\det V_2) + R\pi_{X*}(\pi_1^*\bbO(-1)\otimes \pi_2^*\bbO(3)) \notag \\ &- R\pi_{X*}(\pi_1^*\bbO(-2)\otimes \pi_2^*T\bbP V_2) \notag \\
    &= R\pi_{X*}(\pi_1^*\bbO(-3)\otimes \pi_2^*\det V_2) + R\pi_{X*}(\pi_1^*\bbO(-1)\otimes \pi_2^*\bbO(3))\\
    &- R\pi_{X*}(\pi_1^*\bbO(-2)\otimes \pi_2^*\bbO(1)\otimes V_2) + R\pi_{X*}(\pi_1^*\bbO(-2)). \notag 
\end{align}
We can compute each of the four terms in \eqref{eq_with_4_terms} explicitly to derive a formula in $R(G)$ for $\pi_{X*}e^{K_G}(W)$. The following is well-known, but we include a statement for clarity.

\begin{lemma}\label{lem:cohomology_of_PV}
    Let $V$ be a $G$-representation of dimension $n+1$. Then we have isomorphisms of $G$-representations
    $$
H^i(\mathbb{P}V, \bbO(m))=
\begin{cases}
\Sym^m V^*, \text{ for } i = 0\\
\wedge^{n+1} V, \text{ for } i = n, m=-n-1)\\
\wedge^{n+1} V \otimes \Sym^{-n-1-m} V^*, \text{ for } i = n)\\
0, \text{ for }0<i<n. 
\end{cases}
$$
\end{lemma}

We state another useful lemma: 

\begin{lemma} \label{lem:pshfwd_exceptional}
$R\pi_{2*}\bbO(2E)\cong \bbO_{\bbP V_2}-\bbO_{\bbP V_2}(-6)+\bbO_{\bbP V_2}(-9)^{\oplus 2}-\bbO_{\bbP V_2}(-12)$.
\end{lemma}
\begin{proof}
We have a short exact sequence on X, 
\[
0\to \bbO(2E) \to \bbO(3E) \to i_{E*}\bbO_{E}(-3)\to 0
\]
where $i_E\colon E\hookrightarrow X$ is the inclusion of the exceptional divisor of the blow up of $\bbP V_2$ at $B$, $E = \bbP N_B \bbP V_2$. From this, we see 
\[
R\pi_{2*}i_{E*}\bbO_{E}(-2)\cong i_{B*}R\pi_2|_{E*}\bbO_E(-2)\cong i_{B*}H^1(E, \bbO_E(-2))[1]\cong i_{B*}\wedge^2 N_{B}^*\bbP V_2[1]
\]
where $i_B\colon B\hookrightarrow \bbP V_2$ is the inclusion. Since $N_{B}^*\bbP V_2\cong \bbO(-3)\oplus \bbO(-3)$, we conclude $\wedge^2 N_{B}^*\bbP V_2[1]\cong \bbO(-6)[1]$. 

Let $\mathcal{I}_B$ be the ideal sheaf of $B$. From the short exact sequence 
\[
0\to \bbO(-6)\stackrel{(f,-g)}{\longrightarrow }\bbO(-3)^{\oplus 2}\stackrel{(f,g)}{\longrightarrow}\mathcal{I}_B\to 0
\]
on $X$ it is clear that $i_{B*}\bbO(-6)\cong \bbO(-6)-\bbO(-9)^{\oplus 2}+\bbO(-12)$. Combining these facts, we see 
\[
R\pi_{2*}\bbO(2E)\cong \bbO_{\bbP V_2}-\bbO_{\bbP V_2}(-6)+\bbO_{\bbP V_2}(-9)^{\oplus 2}-\bbO_{\bbP V_2}(-12).
\]
\end{proof}
Note also that $\pi_1^*\bbO(1)\cong \pi_2^*\bbO(3)\otimes \bbO(-E)$. Thus 
\begin{align*}
R\pi_{X*}(\pi_1^*\bbO(-2)) &= R\pi_{X*}(\pi_2^*\bbO(-6)\otimes \bbO(2E) \\ 
&= Rp_{2*}(\bbO(-6)\otimes R\pi_2^*\bbO(2E). 
\end{align*}
By Lemma \ref{lem:pshfwd_exceptional},
\[
R\pi_{X*}(\pi_1^*\bbO(-2)) = Rp_{2*}(\bbO(-6)-\bbO(-12)+\bbO(-15)^{\oplus 2}-\bbO(-18). 
\]
By Lemma \ref{lem:cohomology_of_PV},
\begin{equation}\label{hom_eq:final_form_I}
    R\pi_{X*}(\pi_1^*\bbO(-2)) = \wedge^3V_2\otimes(\Sym^3 V_2^* - \Sym^9 V_2^* + 2\Sym^{12}V_2^* - \Sym^{15}V_2^*). 
\end{equation}

By the same reasoning, 
\begin{align}\label{hom_eq:final_form_II}
R\pi_{X*}(\pi_1^*\bbO(-2) &\otimes \pi_2^*\bbO(1)\otimes V_2) \notag \\
&=V_2\otimes Rp_{2*}(\bbO(1)\otimes R\pi_{2*}(\pi_2^*\bbO(-6)\otimes \bbO(2E))) \notag \\ 
&= V_2\otimes Rp_{2*}(\bbO(-5)-\bbO(-11) +\bbO(-14)^{\oplus 2} -\bbO(-17))\notag \\ 
&= V_2\otimes \wedge^3V_2\otimes (\Sym^2V_2^* - \Sym^8V_2^* +2\Sym^{11} V_2^* - \Sym^{14} V_2^*). 
\end{align}

Since $\pi_1^*\bbO(1)\cong \pi_2^*\bbO(3)\otimes \bbO(-E)$, $R\pi_{X*}(\pi_1^*\bbO(-1)\otimes \pi_2^*\bbO(3)) = R\pi_{X*}\bbO(E)$. Thus 
\begin{equation}\label{hom_eq:final_form_III}
    R\pi_{X*}(\pi_1^*\bbO(-1)\otimes \pi_2^*\bbO(3)) = Rp_{2*}(R\pi_{2*}\bbO(E)) = Rp_{2*}\bbO_{\bbP V_2} = \bbO_{\op{Spec}\C}.
\end{equation}

To compute $R\pi_{X*}(\pi_1^*\bbO(-3)\otimes \pi_2^*\det V_2)$ we use the following lemma, whose proof is omitted as it is similar to the proof of Lemma \ref{lem:pshfwd_exceptional}.

\begin{lemma}\label{lem:pshfwd_O(3E)}
    $R\pi_{2*}\bbO(3E) = R\pi_{2*}\bbO(2E)+\bbO(-12)^{\oplus 4} - \bbO(-9)^{\oplus 2} - \bbO(-15)^{\oplus 2}.$
\end{lemma}

With this is mind, 
\begin{align}\label{hom_eq:final_form_IV}
R\pi_{X*}(\pi_1^*\bbO(-3) &\otimes \pi_2^*\det V_2) \notag \\ &= R\pi_{X*}(\pi_2^*\bbO(-9)\otimes \bbO(3E)\otimes \pi_2^*\det V_2) \notag \\ 
&= \det V_2 \otimes Rp_{2*}(\bbO(-9)\otimes R\pi_{2*}\bbO(3E)) \notag \\ 
&= \det V_2 \otimes Rp_{2*}(\bbO(-9)-\bbO(-15)+\bbO(-21)^{\oplus 3}-\bbO(-24)^{\oplus 2}) \notag \\ 
&= \det V_2\otimes \wedge^3 V_2\otimes (\Sym^6 V_2^*-\Sym^{12}V_2^* + 3\Sym^{18}V_2^* -2\Sym^{21}V_2^*). 
\end{align}

In total we conclude from \eqref{hom_eq:final_form_I}, \eqref{hom_eq:final_form_II}, \eqref{hom_eq:final_form_III}, and \eqref{hom_eq:final_form_IV}, 

\begin{theorem}\label{thm:hom_bundle}
With the assumptions in Theorem \ref{thm:rational_cubic_count}, there is an $R(G)$-valued count of rational cubics through $S$ given by 
\begin{align*}
\pi_Xe^{K_G}(W) &= \sum_{\substack{G\cdot x,\\ d\pi_1(x)=0}} \tr_{G_x}^{G}\deg_x^{G_x}(d\pi_1) \\ 
&= \sum_{\substack{G\cdot x,\\ d\pi_1(x)=0}} \tr_{G_x}^{G}(1) \\
&= \wedge^3V_2\otimes(\Sym^3 V_2^* - \Sym^9 V_2^* + 2\Sym^{12}V_2^* - \Sym^{15}V_2^*) \\ 
&- V_2\otimes \wedge^3V_2\otimes (\Sym^2V_2^* - \Sym^8V_2^* +2\Sym^{11} V_2^* - \Sym^{14} V_2^*) \\ 
&+ \bbO_{\op{Spec}\C}\\
&+ \det V_2\otimes \wedge^3 V_2\otimes (\Sym^6 V_2^*-\Sym^{12}V_2^* + 3\Sym^{18}V_2^* -2\Sym^{21}V_2^*)
\end{align*}
where $W=\op{Hom}(\pi_1^*T^*\bbP V_1, T^*X) \to X$ and $d\pi_1$ is the section of $W$ determined by $\pi_1$.  
\end{theorem}

One checks directly that this has rank 12. By \cite[Lemma 5.21]{brazeuler}, this is independent of choice of section. It's worth noting that while $\pi_{X*}e^{K_G}(W)$ does compute an equivariant count of rational cubics through $S$, it is not equal to $\pi_{X*}e(TX) = \chi^G(X) = N_{3, \CP^2, S}^G$ in $R(G)$, though it differs by a virtual rank 0 class in $R(G)$. 
\newpage

\bibliographystyle{amsalpha}
\bibliography{local_G_degree}

\end{document}